\documentclass{article}
\usepackage{amsfonts}
\usepackage{amsmath}
\usepackage{amssymb}
\usepackage[dvips]{graphicx}
\newtheorem{theorem}{Theorem}
\newtheorem{proposition}[theorem]{Proposition}
\newtheorem{lemma}[theorem]{Lemma}
\newtheorem{corollary}[theorem]{Corollary}
\newtheorem{definition}{Definition}
\newcommand{\tmb}{\tilde{\mathcal{B}}}
\newcommand{\mb}{\mathcal{B}}
\newcommand{\leb}{\mathcal{L}}
\newcommand{\ind}{\mathbf{1}}
\newcommand{\ms}{\mathcal{S}}

\newcommand{\mh}{\mathcal{H}}
\newcommand{\wmh}{\widetilde{\mh}}
\newcommand{\diam}{\textrm{diam}}

\newcommand{\dist}{\textrm{dist}}
\author{Pablo Shmerkin}
\title{Overlapping self-affine sets}
\begin{document}
\maketitle

\begin{abstract}
We study families of possibly overlapping self-affine sets. Our main example is a family that can be considered the self-affine version of Bernoulli convolutions and was studied, in the non-overlapping case, by F.Przytycki and M.Urba{\'n}ski \cite{prurb}. We extend their results to the overlapping region and also consider some extensions and generalizations.
\end{abstract}

\section{Introduction}

A compact, nonempty set $K\subset\mathbb{R}^n$ is called
``Self-affine'' if it is the attractor of an iterated function system of affine maps.
Self-affine sets represent a natural class of
fractal sets. On one hand, they are a natural generalization of
self-similar sets. On the other hand, they
appear in other areas, like theory of tilings and
dynamical systems. They also represent a prototype for attractors
of general, (smooth) non-linear i.f.s. Despite these facts, they remain
rather mysterious, and the study of their dimensional and
topological properties is fraught with difficulties.

One of the most important results available is Falconer's
Theorem from 1988 \cite{falconersa}. It states that for every collection
$T_1,\ldots,T_k$ of linear endomorphisms of $\mathbb{R}^n$ such
that $\|T_i\|<1/2$ for all $1\le i\le k$, there exists a number
$d=d(T_1,\ldots,T_k)$ -called the ``Falconer dimension'' of $\{T_1,\ldots,T_k\}$- such that for almost every
$v_1,\ldots,v_k\in\mathbb{R}^n$ (in the sense of $nk$
dimensional Lebesgue measure), the attractor of the i.f.s.
$\{T_1+v_1,\ldots,T_k+v_k\}$ has Hausdorff and box-counting dimensions
equal to $d$. There is an explicit, albeit difficult to compute, formula for $d$. (Falconer's original version goes with $1/3$ as the bound for the norms; later Solomyak \cite{measanddim} pointed out that $1/2$ works too).

The hypothesis on the norms can be somewhat relaxed, but is essential. The simplest counterexample comes
from the following family of self-affine sets, studied by F.Przytycki and M.Urba{\'n}ski \cite{prurb}:
for $0<\gamma<\lambda<1$, let $T_{\gamma,\lambda}$ be
the linear map given by
\[
T_{\gamma,\lambda}\left(%
\begin{array}{c}
  x \\
  y \\
\end{array}%
\right) = \left(%
\begin{array}{cc}
  \gamma & 0 \\
  0 & \lambda \\
 \end{array}%
\right) \left(%
\begin{array}{c}
  x \\
  y \\
\end{array}%
\right).
\]
Let $K_{\gamma,\lambda}$ be the attractor
of the i.f.s. $\{T_{\gamma,\lambda}-(1,1),T_{\gamma,\lambda}+(1,1)\}$. When $\lambda<1/2$ these sets are easy to analyze,
but the situation becomes much more complicated when $\lambda>1/2$. Recall that the Bernoulli Convolution $\nu_\lambda$ is defined
as the distribution measure of the random sum $\sum_{i=1}^\infty \pm \lambda^i$,
where the signs are chosen indepently with probabilty $1/2$; see \cite{sixtyyears}, \cite{notesbernoulli} for further information on Bernoulli Convolutions.
For us, the main feature of Bernoulli Convolutions is the following theorem of Solomyak \cite{bernoulli}: for almost all $\lambda \in (1/2,1)$, $\nu_\lambda$ is an abolustely continuous measure with an $L^2$ density.

The result of Przytycki and Urba{\'n}ski is the following: Assume that $K_{\gamma,\lambda}$ is totally disconnected, and $1/2<\lambda$.
If $\nu_\lambda$ has Hausdorff dimension $1$ (which, by Solomyak's Theorem, is the case for almost every $\lambda$) then
\[
\dim_H(K_{\gamma,\lambda}) = \dim_B(K_{\gamma,\lambda}) = 1 + \frac{\log(2\lambda)}{\log(1/\gamma)}.
\]
However, if $1/\lambda$ is a Pisot number (i.e., an algebraic number greater than $1$ all of whose algebraic conjugates have modulus less
than $1$), then $\dim_H(K_{\gamma,\lambda})<\dim_B(K_{\gamma,\lambda})$. It is well known that the set of Pisot number
accumulates to $2$; thus this implies that the norm bound in Falconer's Theorem is sharp.

The result of Przytycki and Urba{\'n}ski suggests the following question: what can we say about the dimension and topological properties of $K_{\gamma,\lambda}$ in general; i.e. allowing $K_{\gamma,\lambda}$ to be connected? We will show that the same formula for the dimension holds for almost every $(\gamma,\lambda)$ in the natural region. More precisely, we have:

\begin{theorem} \label{th:main}
For almost $(\gamma,\lambda)$ such that $\gamma\lambda<1/2<\gamma$,
\[
\dim_H(K_{\gamma,\lambda}) = \dim_B(K_{\gamma,\lambda}) = D(\gamma,\lambda),
\]
where $D(\gamma,\lambda)= 1 + \log(2\lambda)/\log(1/\gamma)$.
\end{theorem}

Note that a direct connection between absolute continuity of $\nu_{\lambda}$ and dimension of $K_{\gamma,\lambda}$ is lost. This is natural since for countably many values of $(\gamma,\lambda)$ there is an exact coincidence of cylinders which produces a dimension drop. The condition $\gamma\lambda<1/2$
is also a natural one; for $\gamma\lambda>1/2$ one would expect $K_{\gamma,\lambda}$ to have positive Lebesgue measure, and even non-empty interior. Unfortunately, since transversality holds only in a small region inside $\{\gamma\lambda>1/2\}$, our results here are rather limited.

It is convenient to state the result in terms of measures. Let $\mu_{\gamma,\lambda}$ be the natural self-affine measure supported on $K_{\gamma,\lambda}$; it can be defined in several ways, for instance as the distribution of the random sum
\begin{equation} \label{eq:randomsum}
\mu_{\gamma,\lambda} \sim \sum_{i=0}^\infty  \pm (\gamma^i,\lambda^i),
\end{equation}
where signs are chosen independently with probability $1/2$. Note the close analogy with Bernoulli convolutions; we think of $\mu_{\gamma,\lambda}$ as \textit{self-affine} Bernoulli convolutions.
\begin{theorem} \label{th:mainbis}
There is an open set
\[
\mathcal{U}\subset\{(\gamma,\lambda):0<\gamma<\lambda, \gamma\lambda>1/2\},
\]
containing a neighborhood of $(1,1)$ and of the curve $\{\gamma\lambda=1/2\}$ such that
\begin{enumerate}
\item For almost all $(\gamma,\lambda)\in \mathcal{U}$, $\mu_{\gamma,\lambda}$ is absolutely continuous with an $L^2$ density. In particular, $\leb_2(K_{\gamma,\lambda})>0$ (we will denote $n$-dimensional Lebesgue measure by $\leb_n$).
\item For almost all $(\gamma,\lambda)$ such that $(\gamma^k,\lambda^k)\in \mathcal{U}$ for some $k\ge 2$, $\mu_{\gamma,\lambda}$ is absolutely continuous with a continuous density. In particular, $K_{\gamma,\lambda}$ has nonempty interior.
\end{enumerate}
\end{theorem}

See Corollary \ref{coro:abscontregion}  for the precise definition of $\mathcal{U}$. We remark that our results apply to more general families of self-affine sets, although the theorems above illustrate our main motivation and example. See Theorems \ref{th:generalth} and \ref{th:generalthbis} for the general versions, as well as the extensions and further generalizations presented in Section \ref{sec:extensions}. We stress, however, that we are considering only families of self-affine sets where all the defining maps share the same linear part, which is a diagonalizable map.

The method used to prove Theorems \ref{th:main} and \ref{th:mainbis} is based on the transversality ideas that were successfully applied to many families of self-similar sets, starting with \cite{deleteddigits}; however, some new ideas are needed as well. In particular, we emphasize that the powerful projection scheme developed in \cite{dimexceptions} does not seem to apply in this context, for two reasons. First, we have to deal with two different H\"older exponents simultaneously; second, and more important, the standard notion of transversality does not hold in a large enough region. To overcome the second problem we use transversality concurrently with absolute continuity of Bernoulli convolutions, rather than transversality alone.

Finally, in section \ref{sec:extensions} we consider some additional questions suggested by what is known in the self-similar case. We study some families of exceptions to the almost-everywhere results; most of them are closely related to Pisot numbers. This is to be expected, since for classical Bernoulli convolutions reciprocals of Pisot numbers are the only known parameters that yield a singular measure.

According to Theorem \ref{th:main}, overlaps do not produce a dimension drop in the region $\{\gamma\lambda<1/2\}$ except for a set of zero measure. We show that they do produce a measure drop (i.e. the Hausdorff measure in the critical dimension is $0$) in a big chunk of the overlapping region. This phenomenon was first observed for families of self-similar sets.

It is well known that if a measure is the attractor of a general i.f.s. then it is either singular or absolutely continuous with respect to Lebesgue measure. A more delicate question is whether, when absolutely continuous, it is actually equivalent to Lebesgue measure on the attractor set. This is known to be true for self-similar sets \cite{sixtyyears}, and here we show that the proof in the self-similar case can be adapted to cover many self-affine measures including the natural measures on $K_{\gamma,\lambda}$.

\section{Hausdorff Dimension}

In this section we obtain our main results on Hausdorff and box-counting dimension of certain certain families of self-affine sets, of which the class discussed in the introduction is a particular example.

Let $D=\{d_1,\ldots,d_m\}$ be a set of real numbers which we call ``digits''. We will normalize $D$ so that $0=d_1<d_2<\ldots<d_m$. Associated to $D$ is a family of linear self-similar sets ${K_\lambda}$ for $0<\lambda<1$, where $K_\lambda$ is the attractor of the i.f.s. $\{\lambda x + d_i\}_{i=1}^m$. Furthermore, let $\nu_\lambda$ be the natural self-similar measure supported on $K_\lambda$; i.e. the probability measure defined by the relation
\[
\nu_\lambda = \sum_{i=1}^m \frac{1}{m} (\nu_\lambda\circ\psi_i^{-1}),
\]
where $\psi_i(x) = \lambda x+d_i$.

Now consider a two-dimensional version of this construction, where the i.f.s. is
\[
\left\{  \phi_i \left(%
\begin{array}{c}
  x \\
  y \\
\end{array}%
\right) = \left(%
\begin{array}{cc}
  \gamma & 0 \\
  0 & \lambda \\
 \end{array}%
\right) \left(%
\begin{array}{c}
  x \\
  y \\
\end{array}%
\right) + \left(%
\begin{array}{c}
  d_i \\
  d_i \\
\end{array}%
\right) : 1\le i\le m
\right\}.
\]
Let $K_{\gamma,\lambda}$ be the attractor set, and $\mu_{\gamma,\lambda}$ the natural self-affine measure supported on it. These will be our main object of study in this paper. See Figure 1 for some examples.

\begin{figure}
\begin{center}
\includegraphics[width=0.9\textwidth]{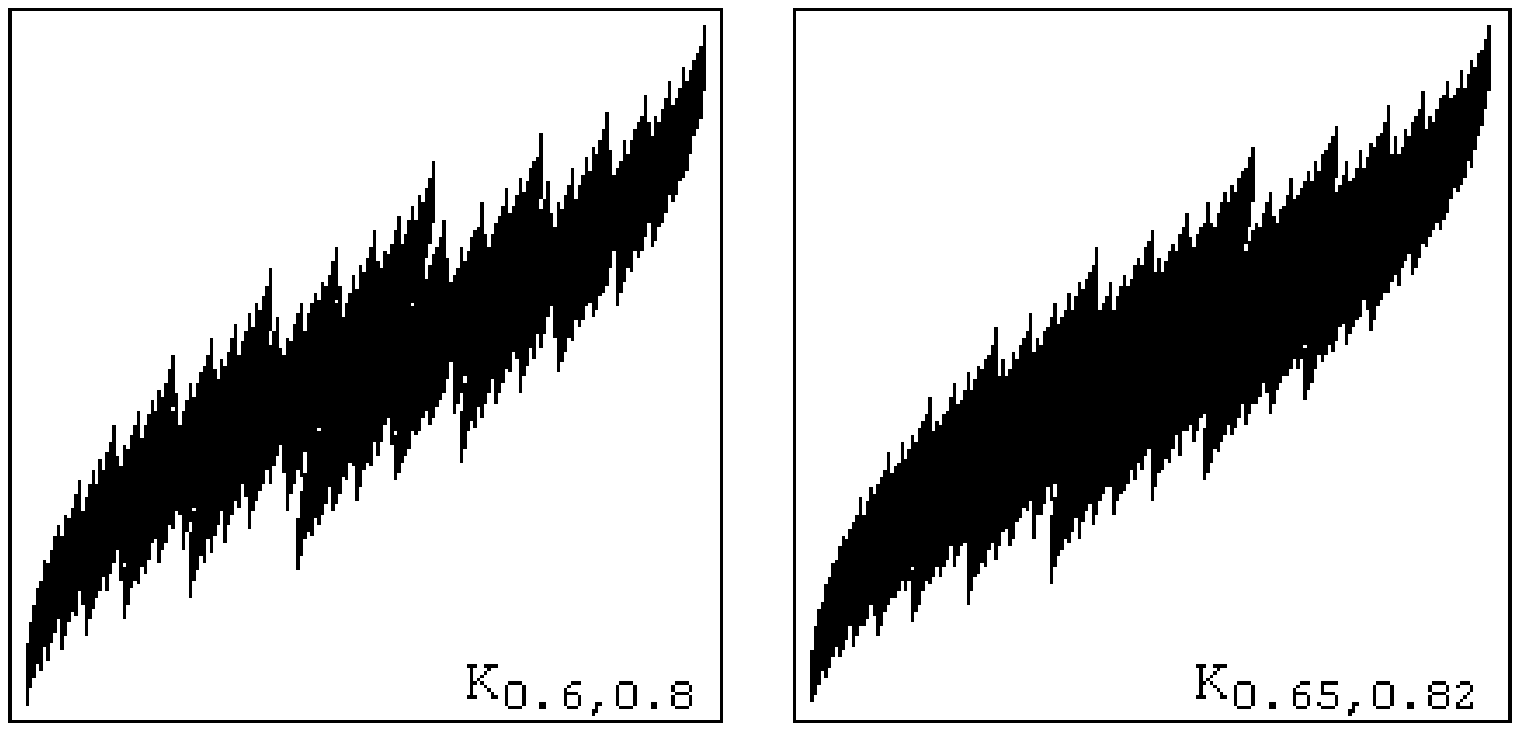}
\end{center}
\textbf{Figure 1}. Both figures correspond to the digit set $\{0,1\}$. Although both pictures look similar, the one on the left actually corresponds to an attractor of dimension strictly less than $2$, while for parameters close to the one on the right we know that the attractor typically has positive Lebesgue measure. In general, self-affine sets are harder to visualize than self-similar sets due to the fact that, after some iterations, cylinders sets have a very large excentricity.
\end{figure}

Our most important example is the set of digits $D=\{0,1\}$. In this case, $\{\nu_\lambda\}$ and $\{\mu_{\gamma,\lambda}\}$ are the classical and self-affine Bernoulli convolutions respectively. More precisely, the measure defined by (\ref{eq:randomsum}) corresponds to the digit set $\{-1,1\}$, but the digits $\{0,1\}$ yield the same measures up to rescaling and translating. For the sake of simplicity we will use the digits $\{0,1\}$ at all places except where we are dealing with the Fourier transform, for which the digits $\{-1,1\}$ yield a slightly simpler formula.

Next we define a relevant class of power series. Let
\[
\tmb = \tmb(D) = \left\{ \sum_{i=0}^\infty c_i \,x^i,\quad c_i \in D-D \right\}.
\]
Let also $\mb$ be the subset of $\tmb$ of power series with non-zero constant term. Moreover, define $\mb_k$ as $x^k\mb$; in other words, $\mb_k$ is the set of power series whose first non-zero coefficient is $c_k$. We will identify (subsets of) $\tmb$ with (subsets of) the symbolic space $(D-D)^\mathbb{N}$.

As a matter of notational convenience we will often write $\Pi_\lambda(f)$ instead of $f(\lambda)$, and $\Pi_{\gamma,\lambda}(f)$ instead of $(f(\gamma),f(\lambda))$. When we do so, we will use Greek letters such as $\omega$ instead of $f$.

We will need a notion of transversality for power series. Our definition is very close to the standard one; however, we need to consider subsets of the class $\mb$ of power series rather than all of $\mb$.

\begin{definition} \label{def:transv}
Let $\mb'$ be a subset of $\mb$. We say that $J\subset(0,1)$ is a set of transversality for $\mb'$ if there exists a constant $M>0$ such that
\begin{equation} \label{eq:transv}
\leb_1(\{\alpha\in J:|g(\alpha)|<r\}) < M\,r
\end{equation}
for every $g\in\mb'$ and $r>0$.

If a set $\tmb'\subset \tmb$ can be expressed as $\bigcup_{j=0}^\infty x^k\,\mb'$ for some $\mb'\subset\mb$, we say that $J$ is a set of transversality for $\tmb'$ if it is a set of transversality for $\mb'$.
\end{definition}

One easy but important observation is that to prove transversality it is enough to show if $f\in\mb'$ then
\[
|f(\alpha)| + |f'(\alpha)| > c,
\]
for all $\alpha\in J$ and some $c>0$.

Let $\rho$ be the uniform Bernoulli measure on $D^\mathbb{N}$; i.e.
\[
\rho = \left( \frac{1}{m},\ldots,\frac{1}{m} \right)^\mathbb{N}.
\]
The projection map $\Pi_\lambda$ from $D^\mathbb{N}$ to $\mathbb{R}$ maps $\rho$ onto the Bernoulli convolution $\nu_\lambda$. The measure $\rho$ induces a measure $\eta$ on $\tmb=(D-D)^\mathbb{N}$ via the difference map. In other words, $\eta$ is the Bernoulli measure on $\tmb$ where the symbol $c\in D-D$ has weight
\[
W(c) = \frac{|\{(a,b)\in D^2:a-b = c\}|}{m^2}.
\]
In particular, $W(0)=1/m$. Let $\sigma$ be the shift operator on $\tmb$.

Finally, let us define some relevant sets. Let
\begin{eqnarray}
\Delta & = & \{ (\gamma,\lambda)\in (0,1)^2: \gamma<\lambda, 1/m<\lambda\};\nonumber \\
\Gamma(a,b) & = & \{ (\gamma,\lambda)\in \Delta: \gamma\lambda \in (a,b)\}\nonumber
\end{eqnarray}
Throughout the paper we will assume that $\gamma$ and $\lambda$ have the same signs, although actually all the results and proofs are valid regardless of signs (the conditions on $\gamma$ and $\lambda$ have to be replaced by conditions on their absolute values; for example $\gamma\lambda<1/m$ becomes $|\gamma||\lambda|<1/m$).

We are now in a position to state a technical proposition that contains our key estimate; in particular, it is only here that transversality gets used.

\begin{proposition} \label{prop:main}
Let $J_1,J_2$ be closed intervals such that $J_1\times J_2\subset \Gamma(0,1/m)$.
\begin{enumerate}
\item
Fix $0<\xi<\varepsilon$. Suppose that for some $\lambda_0\in J_2$ there is a constant $K>0$ such that
\begin{equation} \label{eq:corrbound}
(\nu_{\lambda_0}\times\nu_{\lambda_0})\{ (x,y):|x-y|\le r \} \le K \,r^{1-\xi},
\end{equation}
for all $r>0$.
Assume also that $J_1$ is an interval of transversality for $\tmb'\subset\tmb$. Let $\gamma_0=\min J_1$, and
\begin{equation} \label{eq:dimformula}
D = 1 + \frac{\log(m\lambda_0)}{\log(1/\gamma_0)}
\end{equation}
Then
\begin{equation} \label{eq:keyestimate}
\int_{\tmb'}\int_{J_1} |\Pi_{\gamma,\lambda_0}(\omega)|^{-D+\varepsilon} d\gamma \,d\eta(\omega) < K C(J_1,J_2) <\infty.
\end{equation}
\item
Analogously, assume that (\ref{eq:corrbound}) holds with $\lambda_0$ replaced by $\gamma_0$ for some $\gamma_0\in J_1$, and that $J_2$ is an interval of transversality for $\tmb'\subset\tmb$. Let $\lambda_0=\min J_2$ and $D$ as in (\ref{eq:dimformula}). Then
\[
\int_{\tmb'}\int_{J_2} |\Pi_{\gamma_0,\lambda}(\omega)|^{-D+\varepsilon} d\lambda \,d\eta(\omega) < K C(J_1,J_2) <\infty.
\]
\end{enumerate}
\end{proposition}

\textit{Proof}. We will prove only the first part of the proposition; the second is just a restatement with $\gamma$ and $\lambda$ interchanged. Let $H(\omega)$ denote the inner integral in (\ref{eq:keyestimate}). Using Fubini's Theorem and performing the change of variables $u=\gamma_0^{-t(D-\varepsilon)}$ we get
\begin{eqnarray}
H(\omega) & = & \int_{0}^\infty \leb_1(\{\gamma\in J_2:|\Pi_{\gamma,\lambda_0}(\omega)|^{-D+\varepsilon}>u\})du \nonumber\\
\label{eq:tech1}
& = & \log\gamma_0^{\varepsilon-D}\int_0^\infty \leb_1(\{\gamma\in J_1:|\Pi_{\gamma,\lambda_0}(\omega)|<\gamma_0^t\})\left(\frac{m\lambda_0}{\gamma_0}\right)^t \gamma_0^{\varepsilon t},
\end{eqnarray}
where we used the identity $\gamma_0^{-D} = m\lambda_0/\gamma_0$.

Suppose that $\omega\in\mb_k\cap\tmb'$, and observe that
\[
\Pi_{\gamma,\lambda_0}(\omega) = (\gamma^k\Pi_{\gamma}(\sigma^k\omega),\lambda_0^k \Pi_{\lambda_0}(\sigma^k\omega)),
\]
whence
\[
|\Pi_{\gamma,\lambda_0}(\omega)| <\gamma_0^t \quad\Longrightarrow\quad |\Pi_{\gamma}(\sigma^k\omega)|<\gamma_0^{t-k} \textrm{ and } |\Pi_{\lambda_0}(\sigma^k\omega)|<\gamma_0^t \lambda_0^{-k}.
\]
Therefore we obtain from (\ref{eq:tech1}) that
\begin{eqnarray}
H(\omega) & \le & \log\gamma_0^{\varepsilon-D} \int_0^\infty \ind_{(0,\gamma_0^t \lambda_0^{-k})}(|\Pi_{\lambda_0}(\sigma^k\omega)|)\nonumber\\
\label{eq:tech2}
& & \quad \leb_1(\{\gamma\in J_1: |\Pi_{\gamma}(\sigma^k\omega)| <\gamma_0^{t-k}\}) \left(\frac{m\lambda_0}{\gamma_0}\right)^t \gamma_0^{\varepsilon t} dt.
\end{eqnarray}

Write $G(t,\omega)$ for the integrand in the right hand side of (\ref{eq:tech2}), and define
\begin{eqnarray}
A_k & = & \int_{\mb_k\cap\tmb'} \int_0^k G(t,\omega)dtd\eta(\omega); \nonumber\\
B_k & = & \int_{\mb_k\cap\tmb'} \int_k^\infty G(t,\omega)dtd\eta(\omega). \nonumber
\end{eqnarray}

Note that the integral in (\ref{eq:keyestimate}) is equal to
\[
\log(\gamma_0^{\varepsilon-D}) \sum_{k=0}^\infty A_k + B_k,
\]
so our problem is reduced to estimating the series $\sum_{k} A_k$ and $\sum_k B_k$. Recall that $\nu_{\lambda_0}=\rho\circ\Pi(\cdot,\lambda_0)^{-1}$. Using this, (\ref{eq:corrbound}) and the definition of $\eta$ we see that
\[
\eta(\{\omega\in\tmb:|\Pi_{\lambda_0}(\omega)|<r\}) \le K\,r^{1-\xi} .
\]
Furthermore, since $\eta$ is a product measure,
\begin{equation} \label{eq:tech3}
\eta(\{\omega\in\mb_k:|\Pi_{\lambda_0}(\sigma^k\omega)|<r\}) \le \eta(\mb_k) K r^{1-\xi} = Km^{-k-1}r^{1-\xi}.
\end{equation}
From here and Fubini's theorem we compute
\begin{eqnarray}
A_k & \le & \int_{\mb_k\cap\tmb'}\int_0^k  \ind_{(0,\gamma_0^t \lambda_0^{-k})}(|\Pi_{\lambda_0}(\sigma^k\omega)|)  \left(\frac{m\lambda_0}{\gamma_0}\right)^t \gamma_0^{\varepsilon t} dt d\eta(\omega)\nonumber\\
& = & \int_0^k \left(\frac{m\lambda_0}{\gamma_0}\right)^t \gamma_0^{\varepsilon t} \int_{\mb_k\cap\tmb'}  \ind_{(0,\gamma_0^t \lambda_0^{-k})}(|\Pi_{\lambda_0}(\sigma^k\omega)|) d\eta(\omega)dt  \nonumber\\
& \le & \int_0^k \left(\frac{m\lambda_0}{\gamma_0}\right)^t \gamma_0^{\varepsilon t} K m^{-k-1} (\gamma_0^t\lambda_0^{-k})^{1-\xi}dt \nonumber \\
& = & \frac{K}{m} (m\lambda_0)^{-k} \lambda_0^{\xi k} \int_0^k (m\lambda_0)^t \gamma_0^{(\varepsilon-\xi)t}dt \nonumber\\
& \le & \frac{K}{m\log(m\lambda_0)} \lambda_0^{\xi k}
\le K C_1(J_2) (\min J_2)^{\xi k},
\end{eqnarray}
where in the last step we used that $\gamma_0^{\varepsilon-\xi} \le 1$ and $\lambda_0\ge\min J_2$. Therefore $\sum_{k=0}^\infty A_k< K C'(J_1,J_2)<\infty$.

It remains to show that $\sum_{k=0}^\infty B_k < C''(J_1,J_2)<\infty$. At this point we make use of the transversality hypothesis. Let $\omega\in\mb_k\cap\tmb'$. We use (\ref{eq:transv}) applied to the map $\Pi(\sigma^k\omega,\cdot)$ to estimate
\[
G(t,\omega) \le M\gamma_0^{-k} \ind_{(0,\gamma_0^t\lambda_0^{-k})}(|\Pi_{\lambda_0}(\sigma^k\omega)|)(m\lambda_0)^t \gamma_0^{\varepsilon t}.
\]
Using this, Fubini's theorem and (\ref{eq:tech3}) we obtain
\begin{eqnarray}
B_k & \le & M\gamma_0^{-k} \int_{\mb_k\cap\tmb'} \int_k^\infty \ind_{(0,\gamma_0^t\lambda_0^{-k})}(|\Pi_{\lambda_0}(\sigma^k\omega)|)(m\lambda_0)^t \gamma_0^{\varepsilon t}\,dt\,d\eta(\omega)  \nonumber\\
& = & M\gamma_0^{-k} \int_k^\infty (m\lambda_0)^t\gamma_0^{\varepsilon t} \int_{\mb_k\cap\tmb'} \ind_{(0,\gamma_0^t\lambda_0^{-k})}(|\Pi_{\lambda_0}(\sigma^k\omega)|) d\eta(\omega) dt\nonumber\\
& \le & M K m^{-1} \lambda_0^{\xi k} (m\gamma_0\lambda_0)^{-k} \int_k^\infty (m\gamma_0\lambda_0)^t \gamma_0^{(\varepsilon-\xi)t} dt \nonumber\\
& \le & \frac{MK}{\log(1/(m\gamma_0\lambda_0))} \gamma_0^{(\varepsilon-\xi)k}
< M K C_2(J_1,J_2) (\min J_1)^{(\varepsilon-\xi)k} ,\nonumber
\end{eqnarray}
where we used that $m^{-1}\lambda_0^{\xi k}<1$ and, for the convergence of the last integral, that $m\gamma_0\lambda_0<1$. This completes the proof. $\blacksquare$

\medskip

In order to apply the previous lemma to obtain information about Hausdorff dimension we need to establish transversality. The following easy lemma reduces this problem to the estimation of roots of power series.

\begin{lemma} \label{lemma:transv}
Let $J_1,J_2 \subset (1/m,1)$ be closed intervals such that the following holds: if $\gamma\in J_1,\lambda\in J_2$ and $f\in\mb$, then $\gamma$ and $\lambda$ are not both double roots of $f$. Then $\mb$ can be partitioned into two disjoint subsets $\mb',\mb''$ such that $J_1, J_2$ are sets of transversality for $\mb', \mb''$ respectively.
\end{lemma}
\textit{Proof}. Since $\mb$ is a normal family and $J_1,J_2$ are closed, the following number is well defined and, by hypothesis, positive:
\[
c = \min\{ \max\{|f(\gamma)|,|f'(\gamma)|,|f(\lambda)|,|f'(\lambda)|\}:\gamma\in J_1,\lambda\in J_2,f\in\mb \}.
\]
Let
\[
\mb' = \{f\in\mb: \max\{|f(\gamma)|,|f'(\gamma)|\} \ge c   \}.
\]
Let also $\mb'' = \mb\backslash\mb'$, and note that
\[
\mb'' \subset \{f\in\mb: \max\{|f(\lambda)|,|f'(\lambda)|\} \ge c   \}.
\]
The lemma is now clear from the definition of transversality.
$\blacksquare$

\medskip

Of course, in order to effectively use the previous lemma we need some information on the location of the roots of the power series in $\mb$. This is provided by the following result, which is a simple modification of Theorem 2 in \cite{multroots} (or rather the more general version stated after Theorem 4).
\begin{theorem} \label{th:transvest}
Let
\begin{equation} \label{eq:defchi}
\chi =  \frac{\max_{i,j} |d_i-d_j|}{\min_{i\ne j} |d_i-d_j| }.
\end{equation}
If $0<\alpha_1 \le\ldots\le\alpha_k$ are roots of $f\in\mb$, counted with multiplicity, then
\begin{equation} \label{eq:transvest1}
\prod_{i=1}^k \alpha_i \ge \left(1+\frac{1}{k}\right)^{-k/2} (\chi^2 k+1)^{-1/2}.
\end{equation}
In particular, if
\begin{equation} \label{eq:transvest2}
J_1\times J_2 \subset \left\{ (\gamma,\lambda):\gamma\lambda < \frac{4}{5}(4\chi^2+1)^{-1/4}\right\},
\end{equation}
then the hypothesis of Lemma \ref{lemma:transv} holds.
\end{theorem}
\textit{Proof}. The core of the proof follows closely the proof of Theorem 2 in \cite{multroots}. Let $g\in\mathcal{B}$, and let $f(x) = g(x)/g(0)$; note that $g$ is a monic power series with coefficients bounded by $\chi$. Let
\[
F(x) =  1+\sum_{i=1}^\infty \chi x^{i} = 1+\frac{\chi x}{1-x}.
\]
A straightforward calculation shows that
$\|f(Rz)\|_2^2 \le F(R^2)$, where $\|\cdot\|_2$ denotes $L^2$ norm on the unit circle. Fix $0<R<1$ and Let $k_0 = \max\{j:\alpha_j<R\}$. Then, using Jensen's formula and Jensen's inequality we obtain
\begin{eqnarray}
\sum_{i=1}^{k_0} \log(R/\log\alpha_i) & = & \frac{1}{2\pi} \int_{|z|=1} \log|f(Rz)|dz \nonumber\\
& \le & \log\left(\frac{1}{2\pi} \int_{|z|=1} |f(Rz)|^2dz \right)^{1/2}\nonumber\\
& \le & \log(F(R^2))^{1/2}\nonumber.
\end{eqnarray}
Therefore
\[
\prod_{i=1}^k \alpha_i \ge R^{k-k_0} \prod_{i=1}^{k_0} \alpha_i \ge R^k (F(R^2))^{-1/2}.
\]
Taking $R=(k/(k+1))^{1/2})$ yields (\ref{eq:transvest1}), while setting $k=4, \alpha_1=\alpha_2=\gamma$ and $\alpha_3=\alpha_4=\lambda$ immediately gives (\ref{eq:transvest2}).
$\blacksquare$

\medskip

We will now state the main result of this section. Theorem \ref{th:main} (in fact a more general version) will be obtained as a corollary. We start by recalling the definition of \textit{lower correlation dimension} of a measure $\mu$ on $\mathbb{R}^n$:
\[
\underline{\dim}_2(\mu) = \liminf_{r\rightarrow 0} \frac{\log((\mu\times\mu)\{(x,y):|x-y| \le r \})}{\log r}.
\]
If the above limit exists, we say that the correlation dimension exists and is given by the limiting value.

\begin{theorem} \label{th:generalth}
Assume that:
\begin{enumerate}
\item
There is an open interval $I\subset (1/m,1)$ such that $\nu_\lambda$ has lower correlation dimension $1$ for almost every $\lambda\in I$;
\item
There is an open interval $J\subset (0,1)$ such that if $f\in\mb$, $f$ has no double roots on $J$;
\item
There exists an open region $R$ such that if $(\gamma,\lambda)\in R$, $f\in\mb$, then $\gamma$ and $\lambda$ are not both double roots of $f$.
\end{enumerate}
Let
\begin{eqnarray}
\ms_1 & = & R \cap (I\times I) \cap \Gamma(0,1/m);  \nonumber\\
\ms_2 & = & (J \times I) \cap \Gamma(0,1/m); \nonumber\\
\ms & = & \ms_1 \cup \ms_2\nonumber.
\end{eqnarray}

Then for almost all $(\gamma,\lambda) \in \ms$,
\begin{equation} \label{eq:dimform}
\dim_H(K_{\gamma,\lambda}) = \dim_B(K_{\gamma,\lambda}) = 1 + \frac{\log(m\lambda)}{\log(1/\gamma)}.
\end{equation}
\end{theorem}

Before proving this theorem, some remarks are in order:
\begin{enumerate}
\item Correlation dimension is known to exist for arbitrary self-similar measures, see \cite{lqdimensions}.
\item One can only hope for this theorem to be valid in $\Gamma(0,1/m)$.
Outside this region Falconer's dimension has another expression and one would expect Falconer's dimension to coincide with Hausdorff dimension for almost every parameter. Therefore, one would like to make $I, J$ and $R$ as large as possible; we will see that the case were the digits are equally spaced we do have $\ms=\Gamma(0,1/m)$ (note that for this to hold it is necessary that $\sup I=1$).
\item It is well known that in the self-similar case, if $J$ is as the statement of the theorem, then $\nu_\lambda$ is absolutely continuous for almost all $\lambda\in J\cap \Gamma(1/m,1)$ (see \cite{notesbernoulli}, Theorem 4.3 for a proof); in particular, $J\cap (1/m,1)\subset I$. However, in order to obtain non-trivial results, $I$ has to be larger than $J$. In the case of Bernoulli convolutions, for example, Solomyak's theorem implies that $I=(1/2,1)$ while $J$ is a much smaller interval.
\end{enumerate}

\textit{Proof of Theorem \ref{th:generalth}}. Let us denote the right hand side of (\ref{eq:dimform}) by $D(\gamma,\lambda)$. The inequality $\dim_B(K_{\gamma,\lambda})\le D(\gamma,\lambda)$ is standard for \textit{all} $(\gamma,\lambda)\in\Gamma(1/2,1)$ (this is a particular case of the upper bound in Falconer's theorem). Therefore it is enough to show that for every $\varepsilon>0$ and for all $(\gamma_0,\lambda_0)\in S$ there are closed (non-degenerate) intervals $J_1 \ni \gamma_0, J_2 \ni \lambda_0$ such that
\begin{equation} \label{eq:generalth1}
\dim_H(K_{\gamma,\lambda}) > D(\gamma,\lambda) - 2\varepsilon
\end{equation}
for almost all $(\gamma,\lambda) \in J_1\times J_2$.

\smallskip

We henceforth fix $\varepsilon>0$, $(\gamma_0,\lambda_0)\in \ms$ and choose $J_1,J_2\subset\mathbb{R}$ such that:
\begin{itemize}
\item[$i)$] $\min J_1=\gamma_0$, $\min J_2=\lambda_0$.
\item[$ii)$] $|D(\gamma,\lambda)-D(\gamma',\lambda')|<\varepsilon$ for every $(\gamma,\lambda),(\gamma',\lambda')\in J_1\times J_2$.
\item[$iii)$] $J_1\times J_2 \subset \ms_1$ or $J_1\times J_2\subset \ms_2$.
\end{itemize}

If $J_1\times J_2\subset \ms_1$, let $\{\tmb',\tmb''\}$ be the partition of $\tmb$ given by Lemma \ref{lemma:transv}. Otherwise, let $\tmb'=\tmb,\, \tmb'' = \varnothing$. Note that $\tmb$ is a set of transversality for $J$ (and so is trivially $\varnothing$). This follows from the observation after Definition \ref{def:transv} and the fact that $\tmb$ is a normal family.

Let
\[
J_1^N = \{\gamma\in J_1: (\nu_\gamma\times\nu_\gamma)(\{(x,y):|x-y|\le r\})\le N r^{1-\varepsilon/2} \},
\]
If $J_1\times J_2\subset\ms_1$, define $J_2^N$ analogously; otherwise, let $J_2^N=J_2$. Note that by definition of correlation dimension,
\[
J_i \subset \bigcup_{i=1}^\infty J_i^N, \quad i=1,2.
\]

Therefore it suffices to verify (\ref{eq:generalth1}) for almost every $(\gamma,\lambda)\in J_1^N\times J_2^N$ and then let $N\rightarrow\infty$. We fix $N\in\mathbb{N}$ for the rest of the proof.

We recall the well-known Frostman's Lemma: if there exists a Radon measure $\mu$ supported on a compact set $E\subset\mathbb{R}^n$ such that
\[
\int\int |x-y|^{-s} d\mu(x)d\mu(y) < \infty,
\]
then $\dim_H(E)\ge s$. Recall that $\rho$ is the uniform Bernoulli measure on $D^\mathbb{N}$ and $\mu_{\gamma,\lambda}$ the natural self-affine measure on $K_{\gamma,\lambda}$; it is easy to verify that $\mu_{\gamma,\lambda} = \rho\circ\Pi_{\gamma,\lambda}^{-1}$. Let
\begin{equation} \label{eq:generalth2}
\mathcal{I} = \int_{J_1^N}\int_{J_2^N} \int\int |x-y|^{-D(\gamma,\lambda)+2\varepsilon} d\mu_{\gamma,\lambda}(x) d\mu_{\gamma,\lambda}(y) d\gamma d\lambda.
\end{equation}
By Frostman's lemma, it is enough to show that $\mathcal{I}<\infty$. Let $\mathcal{J}(\gamma,\lambda)$ be the inner double integral in (\ref{eq:generalth2}). Passing to the symbolic space and recalling the definition of $\eta$ we obtain
\begin{eqnarray}
\mathcal{J}(\gamma,\lambda) & = & \int_{D^\mathbb{N}}\int_{D^\mathbb{N}} |\Pi_{\gamma,\lambda}(\omega_1)-\Pi_{\gamma,\lambda}(\omega_2)|^{-D(\gamma,\lambda)+2\varepsilon} d\rho(\omega_1)d\rho(\omega_2) \nonumber \\
& = & \int_{\tmb} |\Pi_{\gamma,\lambda}(\omega)|^{-D(\gamma,\lambda)+2\varepsilon} d\eta(\omega) \nonumber\\
& \le & \int_{\tmb} |\Pi_{\gamma,\lambda}(\omega)|^{-D(\gamma_0,\lambda_0)+\varepsilon} d\eta(\omega),
\end{eqnarray}
where for the last inequality we used condition $ii)$. Let $D=D(\gamma_0,\lambda_0)$. Use Proposition \ref{prop:main} to get
\begin{eqnarray}
\lambda\in J_2^N\Rightarrow\int_{\tmb'}\int_{J_1^N} |\Pi_{\gamma,\lambda}(\omega)|^{-D+\varepsilon} d\gamma d\eta(\omega)& < & N C(J_1,J_2) < \infty;\label{eq:generalth3a}\\
\gamma\in J_1^N\Rightarrow\int_{\tmb''}\int_{J_2^N} |\Pi_{\gamma,\lambda}(\omega)|^{-D+\varepsilon} d\lambda d\eta(\omega) & < & N C(J_1,J_2) < \infty. \label{eq:generalth3b}
\end{eqnarray}
(Note that in the case $J_1\times J_2\subset \ms_2$ , (\ref{eq:generalth3b}) holds trivially because $\tmb''=\varnothing$).

Integrating (\ref{eq:generalth3a}) over $J_2^N$ and (\ref{eq:generalth3b}) over $J_1^N$ and applying Fubini we get
\[
\int_{\tmb_i} \int_{J_1^N}\int_{J_2^N}  |\Pi_{\gamma,\lambda}(\omega)|^{-D+\varepsilon} d\gamma d\lambda d\eta(\omega) < \infty,
\]
for $i=1,2$. Adding, interchanging the order of integration again and recalling (\ref{eq:generalth2}) we conclude
\[
\mathcal{I} < \int_{J_1^N} \int_{J_2^N} \int_{\tmb} |\Pi_{\gamma,\lambda}(\omega)|^{-D+\varepsilon} d\eta(\omega) d\lambda d\gamma <\infty.
\]
This completes the proof. $\blacksquare$

\begin{corollary}
If $D=\{0,1,\ldots,m-1\}$ ($m\ge 2$) then
\[
\dim_H(K_{\gamma,\lambda})=\dim_B(K_{\gamma,\lambda}) = 1+\frac{\log(m\lambda)}{\log(1/\gamma)}
\]
for almost all $(\gamma,\lambda)\in\Gamma(0,1/m)$.
\end{corollary}
\textit{Proof}. In the case $D=\{0,1,\ldots,m\}$, $\nu_\lambda$ is absolutely continuous for almost all $\lambda\in (1/m,1)$. This was proved by Solomyak \cite{bernoulli} in the case $m=2$  and by Simon and T{\'o}th \cite{hajnal}   in the case $m>2$. Therefore we can take $I=(1/m,1)$ in Theorem \ref{th:generalth}.

It is an elementary exercise to verify that
\[
\frac{1}{m} < \frac{4}{5}(4(m-1)^2+1)^{-1/4}
\]
for all $m\ge 2$. Therefore, Theorem \ref{th:transvest} tells us that $R=\Gamma(1/m,1)$ is an appropriate region in Theorem \ref{th:generalth}, and
\[
\ms_1 = \Gamma(1/m,1) \cap \{(\gamma,\lambda):\gamma>1/m\}.
\]
On the other hand, we can take $J=(0,1/m)$ (if $f\in\mb$ then $f$ has no zeros at all in $J$, let alone double zeros). Since $I=(1/m,1)$ we have
\[\ms_2 = \Gamma(1/m,1) \cap \{(\gamma,\lambda):\gamma<1/m\}.
\]
The corollary is now clear. $\blacksquare$

\medskip

We remark that in the region $\ms_2$ a more precise result can be obtained using the method of Przytycki and Urba{\'n}ski; this is due to the fact that the i.f.s. verifies the strong separation condition in that region. However, our proof has the advantage of being more elementary.

It is not difficult to use the method of proof of absolute continuity of Bernoulli convolutions to obtain intervals $I$ for more general digit sets, although in general it may difficult to show that $I=(1/m,1)$. We remark, however,  that Theorem \ref{th:generalth} actually holds for all sufficiently small perturbations of $\{0,1,\ldots,m\}$.

\section{Positive Lebesgue measure}

In the region $\Gamma(1/m,1)$ the Falconer dimension of $K_{\gamma,\lambda}$ is $2$ and we expect $K_{\gamma,\lambda}$ to have positive Lebesgue measure and even non-empty interior. Our main objective in this section is to obtain results analogous to Theorem \ref{th:generalth} in the this region. However, it will be convenient to state our results in terms of measures supported on $K_{\gamma,\lambda}$ rather than the sets themselves. Moreover, we will need to allow for more general measures than the uniform ones considered so far.

Let $\mathcal{M}$ denote the set of Borel probability measures on $D^\mathbb{N}$. For $\omega,\omega'\in D^{\mathbb{N}}$, let $i(\omega,\omega')$ be the length of the longest common initial subsequence of $\omega$ and $\omega'$. We define the (lower) correlation dimension of $\rho\in\mathcal{M}$ as
\[
\underline{\dim}_2(\rho) = \liminf_{k\rightarrow\infty} \frac{-\log((\rho\times\rho)\{(\omega,\omega'):i(\omega,\omega')=k\})}{k\log m}.
\]
One important class of examples are the Bernoulli measures $\mathbf{p}^\mathbb{N}$, where $\mathbf{p}=(p_1,\ldots,p_m)$ is a probability vector. An inspection of the definitions shows that
\[
\underline{\dim}_2(\mathbf{p}^\mathbb{N}) = \sum_{i=1}^m p_i^2.
\]
Let us fix for the moment $\rho\in\mathcal{M}$. The projection map $\Pi_{\gamma,\lambda}$ induces a family of measures $\mu_{\gamma,\lambda}$ supported on $K_{\gamma,\lambda}$, namely $\mu_{\gamma,\lambda} = \rho\circ\Pi_{\gamma,\lambda}(\cdot)^{-1}$. We will also need to redefine $\eta$ to reflect the fact that $\rho$ is now allowed to be a more general measure:
\[
\eta(\Omega) =  (\rho\times\rho)\{ (\omega,\omega'):\omega-\omega' \in \Omega \}, \quad \Omega\subset \tmb.
\]
We can express the correlation dimension in terms of $\eta$ as follows:
\begin{equation} \label{eq:corrdim}
\underline{\dim}_2(\rho) = \liminf_{k\rightarrow\infty} \frac{-\log \eta(\mb_k)}{k\log m}.
\end{equation}

\medskip

We now state a technical proposition which will play an analogous role to that of Proposition \ref{prop:main}.

\begin{proposition} \label{prop:mainbis}
Assume that $\rho$ is a product measure. Let $J_1,J_2$ be closed intervals such that $J_1\times J_2\subset \Gamma(m^{-\underline{\dim}_2(\rho)},1)$.
\begin{enumerate}
\item
Assume that $J_1$ is an interval of transversality for $\tmb'\subset \tmb$. Then there is a constant $C=C(J_1,J_2)$ such that
\[
\limsup_{r\rightarrow 0} \frac{1}{r^2}\int_{\tmb'} \leb_1\{\gamma\in J_1:|\Pi_{\gamma,\lambda}(\omega)|\le r\}\,d\eta(\omega) \le C \left\|\frac{d\nu_{\lambda}}{dx}\right\|_2^2,
\]
for all $\lambda\in J_2$, where $\|d\nu_\lambda/dx\|_2=\infty$ if $\nu_\lambda$ does not have a density in $L^2$.
\item
Analogously, if $J_2$ is an interval of transversality for $\tmb'\subset\tmb$ then there exists $C=C(J_1,J_2)$ such that
\[
\limsup_{r\rightarrow 0} \frac{1}{r^2}\int_{\tmb'} \leb_1\{\lambda\in J_2:|\Pi_{\gamma,\lambda}(\omega)|\le r\}\,d\eta(\omega) \le C\left\|\frac{d\nu_\gamma}{dx}\right\|_2^2,
\]
for all $\gamma\in J_1$.
\end{enumerate}
\end{proposition}
Before proving the proposition, we state a simple lemma we will need.
\begin{lemma}
Let $\nu$ be a measure on $\mathbb{R}$ with an $L^2$ density, which we denote by $f$. Then
\[
\lim_{r\rightarrow 0} \frac{1}{r}(\nu\times\nu)\{(x,y):|x-y|\le r \} = 2 \|f\|_2^2.
\]
\end{lemma}
\textit{Proof of the lemma}. This is standard but we were not able to find a reference, so a proof is provided for the convenience of the reader. Let
\begin{equation} \label{eq:mainbis1}
A(r) = (\nu\times\nu)\{(x,y):|x-y|\le r \}.
\end{equation}
Let also $f_r = f * \phi_r$, where
\[
\phi_r = \frac{\ind_{(-r,r)}}{2r}.
\]
Since $\{\phi_r\}$ is an approximate identity, $f_r\rightarrow f$ in $L^2$. Therefore
\[
\frac{A(r)}{2r} = \frac{1}{2r} \int \nu(B(x,r))d\nu(x) = \int f_r(x) f(x)dx \rightarrow \|f\|_2^2.
\]
The lemma is proved. $\blacksquare$

\textit{Proof of Proposition \ref{prop:mainbis}}. We will prove only the first part, since the second is just a restatement with the parameters interchanged. Let $\omega\in\tmb'\cap\mb_k$. Let us write
\[
\phi(r,\omega) = \leb_1\{\gamma\in J_1: |\Pi_{\gamma,\lambda}(\omega)|\le r\}.
\]
Set $\Omega=\Pi_\lambda^{-1}(0,r\lambda^{-k})\cap\tmb'$ and $\gamma_0=\inf J_1$. We have that
\begin{eqnarray}
\phi(r,\omega) & \le & \ind_{\Omega}(\sigma^k\omega) \leb_1\{\gamma\in J_1:|\Pi_\gamma(\sigma^k\omega)|\le r\gamma_0^{-k}\} \nonumber\\
& \le & c r \gamma_0^{-k}\ind_{\Omega}(\sigma^k\omega)\nonumber,
\end{eqnarray}
for some $c>0$, since $J_1$ is an interval of transversality. From this and the fact that $\rho$ is a product measure we obtain
\begin{eqnarray}
\int_{\tmb'\cap\mb_k} \phi(r,\omega)d\eta(\omega) & \le & c r \gamma_0^{-k} \eta(\{\omega\in\mb_k:|\Pi_{\lambda}(\sigma^k\omega)|<r\lambda^{-k}\}) \nonumber\\
& = & c r \gamma_0^{-k}  \eta(\mb_k) \eta(\{ \omega\in\mb: |\Pi_\lambda(\omega)|<r\lambda^{-k} \}) \label{eq:mainbis2}
\end{eqnarray}
Note however that
\[
\eta(\{ \omega\in\mb: |\Pi_\lambda(\omega)|<r\lambda^{-k} \}) = (\nu_\lambda\times\nu_\lambda)(\{(x,y):|x-y|<r\lambda^{-k}\}).
\]
From this, (\ref{eq:mainbis1}) and (\ref{eq:mainbis2}) we obtain
\[
\limsup_{r\rightarrow 0} \frac{1}{r^2}\int_{\tmb'\cap\mb_k} \phi(r,\omega)d\eta(\omega) \le  c' \|d\nu_\lambda/dx\|_2^2 (\gamma_0\lambda)^{-k}\eta(\mb_k).
\]
Since $(\min J_1 )(\min J_2)>m^{-\underline{\dim}(\rho)}$, we deduce from (\ref{eq:corrdim}) that
\[
\sum_{k=0}^\infty (\gamma_0\lambda)^{-k} \eta(\mb_k) < c''(J_1,J_2)<\infty.
\]
This completes the proof.
 $\blacksquare$

\medskip

We can now state our first result giving regions where $\mu_{\gamma,\lambda}$ is absolutely continuous and, in consequence, $\leb_2(K_{\gamma,\lambda})>0$. Roughly speaking, this theorem follows from Proposition \ref{prop:mainbis} in the same way that Theorem \ref{th:generalth} follows from Proposition \ref{prop:main}. However, when $D=\{0,1\}$ and $\rho=\{1/2,1/2\}^\mathbb{N}$, the theorem gives only a small region where absolute continuity holds. We later extend this region somewhat (in particular, we show it contains a neighborhood of $(1,1)$), although, as mentioned in the introduction, results are far for complete here.

\begin{theorem} \label{th:generalthbis}
Let $\rho$ be a product measure on $D^{\mathbb{N}}$, and let $\mu_{\gamma,\lambda}=\rho\circ\Pi_{\gamma,\lambda}^{-1}$.

Assume that:
\begin{enumerate}
\item
There is an open interval $I\subset (1/m,1)$ such that
\begin{equation} \label{eq:condnorm}
\int_{I'} \|d\nu_\lambda/dx\|_2^2 d\lambda< \infty
\end{equation}
whenever $I'$ is compactly contained in $I$.
\item
There is an open interval $J\subset (1/m,1)$ such that if $f\in\mb$, $f$ has no double roots on $J$;
\item
There exists an open region $R$ such that if $(\gamma,\lambda)\in R$ and $f\in\mb$, then $\gamma$ and $\lambda$ are not both double roots of $f$.
\end{enumerate}
Let
\begin{eqnarray}
\ms_1 & = & R \cap (I\times I) \cap \Gamma(m^{-\underline{\dim}_2(\rho)},1);  \nonumber\\
\ms_2 & = & (J \times I)\cap \Gamma(m^{-\underline{\dim}_2(\rho)},1); \nonumber\\
\ms & = & \ms_1 \cup \ms_2\nonumber.
\end{eqnarray}
Then for almost all $(\gamma,\lambda)\in \mathcal{S}$, $\mu_{\gamma,\lambda}$ is absolutely continuous with a density in $L^2$. \end{theorem}
We make some remarks before the proof.
\begin{enumerate}
\item Note that  $\underline{\dim}_2(\rho)=1/m$ when $\rho$ is the uniform Bernoulli measure $D^{\mathbb{N}}$, but in general there is a gap between the regions given by Theorems \ref{th:generalth} and \ref{th:generalthbis}, even when $D=\{0,1\}$. If the analogy with Bernoulli convolutions holds, $m^{-\underline{\dim}_2(\rho)}$ is the treshold for $L^2$ density when $\rho$ is not uniform, but $\mu_{\gamma,\lambda}$ could still be absolutely continuous even when $\lambda\gamma<m^{-\underline{\dim}_2(\rho)}$.
\item Of course, if $\mu_{\gamma,\lambda}$ is absolutely continuous then $\leb_2(K_{\gamma,\lambda})>0$.
\item The hypothesis on $I$ may appear too strong, but in practice it does in fact follow from the same proof that shows absolute continuity of $\nu_\lambda$ for a.e.$\lambda\in I$.
\end{enumerate}

\textit{Proof of Theorem \ref{th:generalthbis}}. The proof follows the scheme of Theorem 4.3 in \cite{notesbernoulli}. However, we need to introduce the variants we have used before in the proof of Theorem \ref{th:generalth}.

Let $J_1,J_2$ be closed intervals such that $J_1\times J_2\subset \ms_1$ or $J_1\times J_2\subset \ms_2$. In the first case, let us choose $\tmb',\tmb''$ such that $J_i$ is an interval of transversality for $\mb_i$ ($i=1,2$). Otherwise, let $\tmb'=\tmb,\tmb''=\varnothing$.

The first steps mimic the proof \cite{notesbernoulli}, Theorem 4.3, so we only sketch them. Let $\underline{D}(\mu_{\gamma,\lambda})(x)$ be the lower density of $\mu_{\gamma,\lambda}$ at $x\in\mathbb{R}^2$. If we can show that
\[
\mathcal{I} = \int_{J_1}\int_{J_2} \int_{\mathbb{R}^2} \underline{D}(\mu_{\gamma,\lambda})(x) d\mu_{\gamma,\lambda}(x) d\lambda d\gamma,
\]
then the criterion for absolute continuity in \cite{mattila}, section 2.12, will imply that $\mu_{\gamma,\lambda}$ is absolutely continuous for almost all $(\gamma,\lambda)\in J_1\times J_2$ and, moreover, $d\mu_{\gamma,\lambda}/dx \in L^2$.

We can proceed like in \cite{notesbernoulli}, Theorem 4.3 to estimate
\begin{equation} \label{eq:generalthbis1}
\mathcal{I} \le \liminf_{r\rightarrow 0} \frac{1}{\pi r^2} \int_{\tmb}  \leb_2\{(\gamma,\lambda)\in J_1\times J_2: |\Pi_{\gamma,\lambda}(\omega)|\le r \} \,d\eta(\omega).
\end{equation}
Therefore it is enough to show that the integral in the right hand side above is bounded by $C_2 r^2$ when restricted to $\tmb'$ or $\tmb''$. We consider only the restriction to $\tmb'$ since the other case is similar. Using Fubini, Proposition \ref{prop:mainbis} and (\ref{eq:condnorm}), we can estimate the integral in (\ref{eq:generalthbis1}) restricted to $\mb_1$ as
\begin{eqnarray}
\int_{J_2} \int_{\tmb'}\leb_1\{\gamma\in J_1: |\Pi_{\gamma,\lambda}(\omega)|\le r \} d\lambda d\eta(\omega) & \le &
C r^2 \int_{J_2} \|d\nu_\lambda/dx\|_2^2 d\lambda\nonumber \\
& \le & C_2 r^2 \nonumber.
\end{eqnarray}
The proof is now complete. $\blacksquare$

\medskip

We now give a concrete region where we can show that $\mu_{\gamma,\lambda}$ is absolutely continuous in the case $m=2$. Other values of $m$ can be handled analogously.

\begin{corollary} \label{coro:abscontregion}
Let
\[
\mathcal{S} = \Gamma\left(1/2,4\times 5^{-5/4}\right) \cup \{ (\gamma,\lambda)\in \Gamma(1/2,1): \gamma < 0.649 \}.
\]
Further, let
\begin{eqnarray}
\mathcal{V} & =  & \left\{ (\gamma,\lambda): \left(\gamma^k,\lambda^k\right)\in S \textrm{ for some } k\ge 2\right\};\nonumber\\
\mathcal{U} & = & \mathcal{S} \cup \mathcal{V}\nonumber.
\end{eqnarray}
Then $\mu_{\gamma,\lambda}$ is absolutely continuous for almost every $(\gamma,\lambda)\in \mathcal{U}$. Furthermore, $\mu_{\gamma,\lambda}$ has a continuous density for almost every $(\gamma,\lambda)\in \mathcal{V}$.
\end{corollary}
\textit{Proof}. It was proved in \cite{bernoulli} that we can take $I=(1/2,1-\varepsilon)$ for any $\varepsilon$ and $J=(0,0.649)$. On the other hand, Theorem \ref{th:transvest} applied with $k=4$ and $\chi=1$ shows that we can take $R=\Gamma(1/2,4 \times 5^{-5/4})$. Therefore Theorem \ref{th:mainbis} shows that $\mu_{\gamma,\lambda}$ is absolutely continuous for almost every $(\gamma,\lambda)\in \mathcal{S}$.

For the rest of the corollary it suffices to show that if $\mu_{(\gamma^k,\lambda^k)}$ is absolutely continuous with a density in $L^2$, then $\mu_{\gamma,\lambda}$ is absolutely continuous with a continuous density. This is standard but we include a proof for completeness.

By decomposing the measure $\mu_{\gamma,\lambda}$ as an infinite convolution of Bernoulli measures we obtain
\begin{eqnarray}
\widehat{\mu}(x,y) & = & \prod_{i=0}^\infty \frac{1}{2} \left( \exp(- i (\gamma^i,\lambda^i)\cdot (x,y))+\exp(- i (-\gamma^i,-\lambda^i)\cdot (x,y)) \right) \nonumber\\
&  = & \prod_{i=0}^\infty \cos(-(\gamma^i x+\lambda^i y)) \label{eq:fourier}.
\end{eqnarray}
(Recall that we are using the digits $\{-1,1\}$ rather than $\{0,1\}$. Therefore, for all $k\ge 2$,
\begin{equation} \label{eq:fourierdecomposition}
\widehat{\mu}_{\gamma,\lambda}(x,y) = \prod_{i=0}^{k-1} \widehat{\mu}_{\left(\gamma^k,\lambda^k\right)}(\gamma^i x,\lambda^i y).
\end{equation}
Assume $\widehat{\mu}_{\left(\gamma^k,\lambda^k\right)}$ has an $L^2$ density; then so does its Fourier transform. Now (\ref{eq:fourierdecomposition}) shows that $\widehat{\mu}_{\gamma,\lambda} \in L^{2/k}$; on the other hand, $\widehat{\mu}_{\gamma,\lambda}$ is bounded since $\mu_{\gamma,\lambda}$ is a finite measure. Therefore $\widehat{\mu}_{\gamma,\lambda} \in L^1$ and, after taking inverse Fourier transform, we conclude that $\mu_{\gamma,\lambda}$ is absolutely continuous with a continuous density, as desired $\blacksquare$

The regions $\mathcal{S}, \mathcal{V}$ are pictured in Figure 2.
\begin{figure}
\begin{center}
\includegraphics[width=0.9\textwidth]{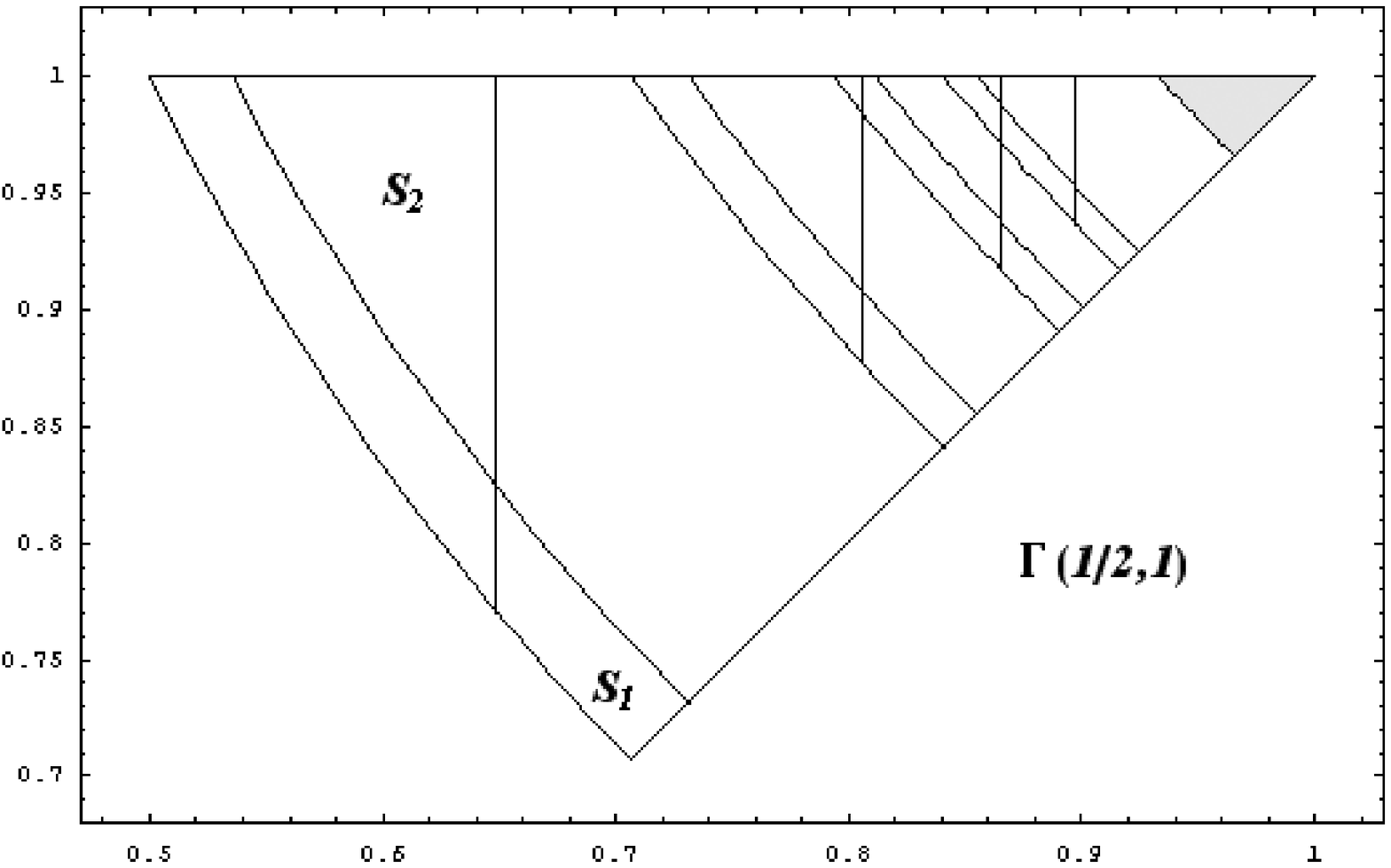}
\end{center}
\textbf{Figure 2}. The region $\mathcal{S}=S_1\cup S_2$ in Corollary \ref{coro:abscontregion}, and the first $3$ pieces or $\mathcal{V}$. The shaded region is $\Gamma(2^{-1/10},1)$, which is also contained in $\mathcal{V}$.
\end{figure}

Several remarks are in order.
\begin{enumerate}
\item $\mathcal{V}$ contains a neighborhood of $(1,1)$; more precisely, $\Gamma(2^{-1/10},1) \subset \mathcal{V}$. Indeed, one can see that
\[
\bigcup_{k=10}^\infty \left((1/2)^{-1/k}, \left(4 \times 5^{-5/4}\right)^{1/k}\right) = (2^{-1/10},1),
\]
as each interval in the union in the left hand side overlaps with the next.
\item Using tricks such as those in \cite{bernoulli} it is possible to extend the region $\mathcal{U}$, but not significantly. It is also possible to use a computer-assisted rigorous estimation of the region $R$ to show that $\mathcal{U}$ can be enlarged to cover more than $90\%$ of $\Gamma(1/2,1)$; details will be discussed in a forthcoming paper.
\end{enumerate}

\section{Extensions and remarks} \label{sec:extensions}

In this section we discuss several natural questions. Some of them can be easily answered with the current techniques, while others appear to be harder. For the most part we restrict ourselves to the case $D=\{0,1\}$ or $D=\{0,1,\ldots,m-1\}$. An exception is Subsection \ref{sec:puretype}, where we consider more general self-affine sets and measures.

We introduce some notation we will be using repeatedly. Recall that
\[
\phi_i(x,y) =  (\gamma x,\lambda y) + (d_i,d_i).
\]
If $u\in D^k$, let $[u]$ be the associated cylinder set in $D^\mathbb{N}$, let $K_{\gamma,\lambda}(u)$ be the cylinder set $\Pi_{\gamma,\lambda}([u])$ and, finally, let
\[
\phi_u = \phi_{u(1)}\circ\ldots\circ \phi_{u(k)}.
\]
Note that $K_{\gamma,\lambda}(u) = \phi_u(K_{\gamma,\lambda})$. We will often omit the subscripts $(\gamma,\lambda)$ whenever they are fixed in a context.

\subsection{Connectedness Loci}

In this subsection we consider the case $D=\{0,1\}$. A natural question is for what values of $(\gamma,\lambda)$ is $K_{\gamma,\lambda}$ connected; this parameter set is called the ``Mandelbrot set'' or ``Connectedness locus'' associated to the family; it will be denoted by $\mathcal{M}$. B. Solomyak \cite{connectloci} has some interesting results on this topic. Without going into details, let us state some of the known basic facts which will be useful later.
\begin{lemma} \label{lemma:conneclocus}
\begin{enumerate}
\item $(\gamma,\lambda)\in\mathcal{M}$ if and only if there is $f\in\mb$ such that $f(\gamma)=f(\lambda)=0$.
\item If $m=1$, $(\gamma,\lambda)\notin\mathcal{M}$ for all $(\gamma,\lambda)$ such that $\lambda<0.649$.
\item $(\gamma,\lambda)\in\mathcal{M}$ for all $(\gamma,\lambda)$ such that $\gamma\lambda>1/2$.
\end{enumerate}
\end{lemma}
\textit{Proof}. \textit{1}. One can easily see that $K_{\gamma,\lambda}$ is connected only if $\Pi_{\gamma,\lambda}$ is not one to one; and when $D=\{0,1\}$ it is if and only if. Since $\tmb = D^\mathbb{N} - D^\mathbb{N}$ the assertion follows.

\textit{2}. This follows immediately from \textit{1}. and the fact that $(0,0.649)$ is an interval of transversality in the strong sense that $f\in\mb$ cannot have two zeros on $(0,0.649)$; see \cite{simplebernoulli}.

\textit{3}. Assume $K=K_{\gamma,\lambda}$ is disconnected. Denote by $A(\delta)$ the $\delta$ neighborhood of $A$. Let $K_0,K_1$ be the cylinders of step $1$, and pick $\delta$ such that $K_0(\delta)$ and $K_1(\delta)$ are disjoint. Since $\phi_i(K(\delta))\subset K_i(\delta)$ because $\phi_i$ is contractive, we must have
\[
\leb_2(K(\delta)) > \leb_2(K_0(\delta)) + \leb_2(K_1(\delta)) = (\det\phi_0+\det\phi_1) \leb_2(K(\delta));
\]
hence $\det\phi_i=\gamma\lambda<1/2$. $\blacksquare$

\subsection{Pisot numbers and the set of exceptions}

For Bernoulli convolutions, the only values of $\lambda$ for which is known that $\nu_\lambda$ is singular are reciprocals of Pisot numbers. This is proved by showing that the Fourier transform $\widehat{\nu}_\lambda(x)$ does not go to zero as $x\rightarrow\infty$. In fact, it is known that $\widehat{\nu}_\lambda$ does not converge to $0$ at infinity if and only if $1/\lambda$ is Pisot (\cite{pisotsalem}, Proposition 15.3.2). The problem of determining all $\lambda$ such that $\nu_\lambda$ is singular is currently open and seems to be very hard.

In \cite{prurb} it is proved that for $\gamma=1/2$ and $\lambda$ the reciprocal of a Pisot number, it is verified that
\[
\dim_H(K_{\gamma,\lambda}) < \dim_B(K_{\gamma,\lambda}) = 1+\frac{\log(2\lambda)}{\log(1/\gamma)}.
\]
(Here $K_{\gamma,\lambda}$ is the attractor for the digit set $\{0,1\}$). Their proof extends readily to arbitrary $\gamma$ such that $K_{\gamma,\lambda}$ is totally disconnected. Although they do not write down the details, they remark that a simpler proof of the same fact can be obtained using the technique of McMullen \cite{mcmullen}.

Here we investigate the set of exceptions to the almost everywhere results obtained in the earlier sections. In particular, we write down the details of the proof suggested by Przytycki and Urba{\'n}ski, which extends to the overlapping case, and in fact it provides examples of $(\gamma,\lambda)\in \Gamma(1/2,1)$ such that $\dim_H(K_{\gamma,\lambda})<2$. Throughout the section we assume that $m=2$, although for the most part analogous considerations are valid for $m>2$ as well.

The first lemma shows that exact coincidence of cylinders produce a dimension drop, and this happens on a dense set.

\begin{lemma}
Let $Q$ be the set of all $(\gamma,\lambda)\in \Gamma(0,1/2)$ such that $\gamma$ and $\lambda$ are roots of a \em polynomial \em with coefficients in $\{-1,0,1\}$. For all $(\gamma,\lambda)\in Q$,
\[
\overline{\dim}_B(K_{\gamma,\lambda}) < 1+\frac{\log(2\lambda)}{\log(1/\gamma)}.
\]
Moreover, the set $Q$ is dense in $\Gamma(0,1/2)\cap\mathcal{M}$, where $\mathcal{M}$ is the connectedness locus for the family $\{K_{\gamma,\lambda}\}$.
\end{lemma}
\textit{Proof}. If $(\gamma,\lambda)\in Q$ we can find two words $u,v\in \{0,1\}^k$ for some $k\ge 1$ such that $\phi_u = \phi_v$ and therefore $K_u = K_v$. It follows that $K$ can be covered by $2^k-1$ rectangles of size $C_1\gamma^k\times C_1\lambda^k$, and more generally, by $(2^k-1)^n$ rectangles of size $C_1\gamma^{nk}\times C_1\lambda^{nk}$. By subdividing each of those rectangles into $C_2 (\lambda/\gamma)^{nk}$ squares of side $\gamma^{nk}$ we conclude that $K_{\gamma,\lambda}$ can be covered by
\[
N(k) = C_3 (2^k-1)^n (\lambda/\gamma)^{nk}
\]
balls of radius $\gamma^{n k}$. Therefore
\[
\overline{\dim}_B(K_{\gamma,\lambda}) \le \limsup_{k\rightarrow\infty} \frac{\log(N(k))}{nk\log(\gamma)} = \frac{\log((2^k-1)\lambda^k)}{k\log(1/\gamma)},
\]
and the first assertion follows. The fact that $Q$ is dense in $\mathcal{M}$ is an immediate consequence of Lemma \ref{lemma:conneclocus} and Rouch\'{e}'s Theorem. $\blacksquare$

We will say that $(\alpha,\beta)$ is a \em Pisot pair \em if $\alpha>1,\beta>1$ and there exists a monic irreducible polynomial $P\in \mathbb{Z}[x]$ such that $\alpha$ and $\beta$ are the only roots of $P$ with modulus greater or equal than $1$. Pisot pairs (or rather more general ``Pisot families'') have been studied by several authors; see for instance \cite{pisotfamily}.

\begin{lemma}
If either $\gamma^{-1}$ or $\lambda^{-1}$ is a Pisot number, or if $(\gamma^{-1},\lambda^{-1})$ is a Pisot pair, then $\widehat{\mu}_{\gamma,\lambda}(x,y)$ does not converge to $0$ as $(x,y)\rightarrow\infty$. In particular, $\mu_{\gamma,\lambda}$ is singular.
\end{lemma}
\textit{Proof}. We know that the vertical and horizontal projections of $\mu_{\gamma,\lambda}$ are $\nu_{\lambda}$ and $\nu_{\gamma}$ respectively. Therefore the restrictions of $\widehat{\mu}_{\gamma,\lambda}$ are $\widehat{\nu}_{\lambda}$ and $\widehat{\nu}_{\gamma}$. Hence whenever one of $\widehat{\nu}_{\lambda}$ or  $\widehat{\nu}_{\gamma}$ does not converge to $0$ at infinity, the same happens to $\widehat{\mu}_{\gamma,\lambda}$. In particular, this is the case if either $\gamma^{-1}$ or $\lambda^{-1}$ are Pisot.

If $(\gamma^{-1},\lambda^{-1})$ is a Pisot pair then we can apply the technique used to prove that $\widehat{\nu}_\lambda(x)$ does not converge to $0$ at infinity. An even closer example is the family of complex Bernoulli convolutions studied in \cite{mandbernoulli}. We indicate the idea, but refer to the proof of Theorem 2.3 in \cite{mandbernoulli} for the details.

Let $P$ be the polynomial in the definition of Pisot pair, and let $\alpha=\gamma^{-1},\beta=\lambda^{-1},\zeta_1,\ldots,\zeta_j$ be the roots of $P$ (so $|\zeta_i|<1$). Since for all integers $n$,
\[
\alpha^n + \beta^n + \sum_{i=1}^j \zeta_j^n \in \mathbb{Z},
\]
we have that $\dist(\alpha^n+\beta^n,\mathbb{Z})< c^n$ for some $c<1$. This fact combined with the expression (\ref{eq:fourier}) for the Fourier transform imply (after some technical considerations) that
\[
|\widehat{\mu}_{\gamma,\lambda}(\pi \gamma^{-N},\pi \lambda^{-N})| > \delta > 0,
\]
for all positive integers $N$ and some $\delta$ independent of $N$. $\blacksquare$

Note that the singularity of $\mu_{\gamma,\lambda}$ in the above lemma is significant only for $(\gamma,\lambda)\in\Gamma(1/2,1)$, for otherwise we know that $\dim_H(K_{\gamma,\lambda})<2$.

It is not difficult to obtain many examples of Pisot pairs using a computer. For instance, the polynomial
\[
P(x) =x^6 - x^5 - x^4 - x^3 + x^2 + x + 1
\]
has exactly two positive roots $\beta>\alpha>1$ and all the other roots are complex and of absolute value less than $1$. The approximate values are $\lambda=\alpha^{-1}=0.754878\ldots$ and $\gamma=\beta^{-1}=0.682328\ldots$. Note that in this example $\gamma\lambda>1/2$.

We finish this section with the extension of the Theorem of Przytycki  and Urba{\'n}ski mentioned in the introduction.
\begin{theorem} \label{th:dimdrop}
Fix $\lambda>1/2$ such that $1/\lambda$ is a Pisot number. There exists a strictly positive continuous function $\theta=\theta_\lambda$ defined on $(0,\lambda)$ such that
\[
\dim_H(K_{\gamma,\lambda}) \le 1+\frac{\log(2\lambda)}{\log(1/\lambda)} - \theta(\gamma).
\]
\end{theorem}
\textit{Proof}. The proof consists of two parts. In the first part we obtain an upper bound for the dimension of $K_{\gamma,\lambda}$ using McMullen's technique; in the second part we show that this upper bound verifies the inequality in the theorem. In the course of the proof, $C_1,C_2$ etc will denote constants independent of $k$ and $\gamma$ (they may depend on the fixed number $\lambda$).

Let
\[
\mathcal{P}_k = \{ \Pi_\lambda(u) : u\in\{0,1\}^k \}.
\]
We recall Garsia's lemma \cite{garsia}: the distance between any two different elements of $\mathcal{P}_k$ is bounded below by $C_1 \lambda^k$. On the other hand, since $\lambda>1/2$ one can easily see, using for instance the greedy algorithm, that the distance between consecutive elements of $\mathcal{P}_k$ is at most $\lambda_k$. Since $0\le \Pi_\lambda(u) \le \lambda/(1-\lambda)$ for all $u\in\{0,1\}^k$, it follows that
\begin{equation} \label{eq:boundNk}
C_2 \lambda^{-k} < \#\mathcal{P}_k < C_3 \lambda^{-k}.
\end{equation}
Let $t_1 < t_2 <\ldots < t_{N_k}$ be the elements of $\mathcal{P}_k$ (so that $\#\mathcal{P}_k=N_k$), and let
\[
a_k(j) = \#\left\{ u \in \{0,1\}^k : \Pi_\lambda(u)=t_j \right\}.
\]
Note that
\begin{equation} \label{eq:sumakj}
\sum_{j=1}^{N_k} a_k(j) = \#\{ 0,1\}^k = 2^k.
\end{equation}
We will show the following:
\begin{equation} \label{eq:mcmullenest}
\dim_H(K_{\gamma,\lambda})) \le \frac{\log\left(\sum_{j=1}^{N_k} a_k(j)^q\right)}{k \log(1/\lambda) },
\end{equation}
where $q=\log\lambda/\log\gamma$. In the course of the proof of (\ref{eq:mcmullenest}) all numbers $\lambda$, $\gamma$ and $k$ will remain fixed.

Let $\Omega$ be the symbolic space $D^\mathbb{N}$, where $D=\{0,1\}^k$. Further let
\[
Z = Z_k(q) = \sum_{j=1}^{N_k} a_k(j)^q.
\]
Let $\omega\in\Omega$. Define the $n$-th symbolic approximate square $S_n(\omega)\subset\Omega$ as $S_n(\omega) = S'_n(\omega)\cap S''_n(\omega)$, where
\begin{eqnarray}
S'_n(\omega) & = & \{ \theta\in\Omega: \theta_i = \omega_i \textrm{ for } 1\le i \le q n  \};\\
S''_n(\omega) & = & \{ \theta\in\Omega: \Pi_\lambda(\theta_i) = \Pi_\lambda(\omega_i) \textrm{ for } q n < i \le n \}.
\end{eqnarray}
Let $\Sigma:\Omega\rightarrow K_{\gamma,\lambda}$ be the projection map given by
\[
\Sigma(\omega) = \left(\sum_{i=1}^\infty \Pi_\gamma(\omega_i)\gamma^{k i},\sum_{i=1}^\infty \Pi_\lambda(\omega_i)\lambda^{k i} \right).
\]
(So in other words $\Sigma(\omega)=\Pi_{\gamma,\lambda}(\omega')$, where $\omega'$ is the sequence obtained by concatenating all the $\omega_i$). Note that $\Sigma$ is surjective but not necessarily injective. Observe that if $\theta\in S_n(\omega)$ then
\begin{eqnarray}
|(\Sigma(\theta) - \Sigma(\omega))_1| & \le &  \sum_{i=\lfloor qn \rfloor+1}^\infty C_4 \gamma^{ni} < C_5 \lambda^{nk},\nonumber\\
|(\Sigma(\theta) - \Sigma(\omega))_2| & \le &  \sum_{i= n}^\infty C_6 \lambda^{ni} = C_7 \lambda^{nk},
\end{eqnarray}
where $(\cdot)_1, (\cdot)_2$ represent the first and second coordinates. Therefore $\Sigma(S_n(w))$ is contained in a ball of center $\Sigma(\omega)$ and radius comparable to $\lambda^{kn}$.

Let $\rho$ be the Bernoulli measure on $\Omega$ giving weight $Z^{-1} a_k(j)^{q-1}$ to all $u$ such that $\Pi_{\lambda}(u) = t_j$ (See \cite{mcmullen} for a motivation for this choice of weights). Our next step is to show that for all $\omega\in\Omega$,
\begin{equation} \label{eq:Snest}
\limsup_{n\rightarrow\infty} \, \frac{\log\rho(S_n(\omega))}{n} + \log Z \ge 0.
\end{equation}
We will do so by using a clever trick due to McMullen. For $u\in D$ write $a(u) = a(j)$ if $\Pi_\lambda(u)=t_j$. Note from the definition of $S_n(\omega)$ that
\[
\log\rho(S_n(\omega)) = \sum_{i=1}^{\lfloor nq \rfloor} \left((q-1)\log(a(\omega_i)) - \log Z\right) + \sum_{i=\lfloor nq \rfloor+1}^{n} \left(q\log(a(\omega_i)) - \log Z\right).
\]
Therefore
\begin{equation} \label{eq:estmm1}
\frac{1}{q n}\left(\log\rho(S_n(\omega)) + n\log Z\right) = \frac{1}{n}\sum_{i=1}^n  \log(a(\omega_i)) - \frac{1}{n q}\sum_{i=1}^{\lfloor nq \rfloor}  \log(a(\omega_i)).
\end{equation}
Write the right hand side above as $S_n/n - S_{qn}/qn$, where
\begin{equation} \label{eq:estmm2}
S_\alpha = \sum_{i=1}^{\lfloor \alpha \rfloor} \log a(\omega_i).
\end{equation}
Note that $S_\alpha/\alpha$ is bounded over all positive $\alpha$. Therefore, by telescoping over the sequence $1,q^{-1},q^{-2}$ we deduce that
\[
\sum_{i=1}^N \left(\frac{S_{q^{-(i+1)}}}{q^{-(i+1)}} - \frac{S_{q^{-i}}}{q^{-i}}\right).
\]
is also bounded over all $N$. Observe that since $S_\alpha-S_{\lfloor \alpha \rfloor}$ is bounded, $S_\alpha/\alpha - S_{\lfloor \alpha \rfloor}/\lfloor \alpha \rfloor \rightarrow 0$ as $\alpha\rightarrow\infty$. Therefore we must have
\[
\limsup_{n\rightarrow\infty} \left(\frac{S_n}{n}-\frac{S_{qn}}{qn}\right) \ge0.
\]
This together with (\ref{eq:estmm1}) and (\ref{eq:estmm2}) show that (\ref{eq:Snest}) is verified.

Recall that $\Sigma(S_n(\omega))$ is contained in a ball $B(\Sigma(\omega),C_8 \lambda^{kn})$. It follows from $(\ref{eq:Snest})$ that if $\mu$ is the projection of $\rho$  under $\Sigma$  then
\[
\liminf_{n\rightarrow\infty} \frac{\log\mu(B(\Sigma(\omega),C_8\lambda^{kn}))}{\log(\lambda^{kn})} \le \frac{\log Z}{-k\log \lambda}.
\]
Since $\Sigma$ is surjective, the mass distribution principle (\cite{falconer3}, Proposition 2.2) shows that (\ref{eq:mcmullenest}) is satisfied. This concludes the first part of the proof.

From now on, $\lambda$ will remain fixed, but we will consider both $k$ and $\gamma$ (or rather $q$) as variables. Let $p_k(j) = 2^{-k} a_k(j)$. Recalling that $q=\log\lambda/\log\gamma$, write the upper bound in (\ref{eq:mcmullenest}) as
\[
U_k(q) = \frac{\log\left(\sum_{i=1}^{N_k} p_k(j)^q\right)}{k\log (1/\lambda)} + \frac{\log 2}{\log(1/\gamma)}.
\]
Let $W_k(q) = \sum_{j=1}^{N_k} p_k(j)^q$. Notice that
\begin{equation} \label{eq:dimdiff}
\liminf_{k\rightarrow\infty} U_k(q) - \left(1+\frac{\log(2\lambda)}{\log(1/\gamma)}\right) = \liminf_{k\rightarrow\infty} \frac{\log W_k(q)}{k\log(1/\lambda)} - (1-q).
\end{equation}
Therefore in order to establish the theorem it is enough to show that the right hand side above defines a continuous function of $q$ which is strictly negative on $(0,1)$.

We claim that for fixed $q$ the sequence $\{W_k(q)\}_k$ is submultiplicative. Indeed,
\[
W_{k_1}(q)W_{k_2}(q) = \sum_{j_1=1}^{N_{k_1}}\sum_{j_2=1}^{N_{k_2}} (p_{k_1}(j_1)p_{k_2}(j_2))^q.
\]
On the other hand, each $p_{k_1+k_2}(j)$ is the sum of one or more numbers of the form $p_{k_1}(j_1)p_{k_2}(j_2)$, and each pair $(j_1,j_2)$ appears in exactly one of the $p_{k_1+k_2}(j)$. The submultiplicativity is then consequence of the inequality
\[
\left(\sum_j r_j\right)^q \le \sum_j r_j^q
\]
for a finite collection of positive numbers $\{r_j\}$, which holds since $0<q<1$. By taking logarithms and using subadditivity we obtain
\[
\lim_{k\rightarrow\infty}  \frac{\log W_k(q)}{k\log(1/\lambda)} = \inf_k \frac{\log W_k(q)}{k\log(1/\lambda)}.
\]
Denote the limiting function by $\tau(q)$. Since $W_k(0)=\#\mathcal{P}_k$ and $W(1)=\sum_j p_k(j)$, we deduce from (\ref{eq:boundNk}) and (\ref{eq:sumakj}) that $\tau(0)=1, \tau(1)=0$.

Also, notice that $\log(W_k(q))$ is a convex function. Since, because of subadditivity, $\tau(q)$ is then a pointwise limit of decreasing convex functions,  $\tau(q)$ must itself be convex. In particular, $\tau$ is continuous and, since it agrees with the linear function $q-1$ at $0$ and $1$, we must have $\tau(q) \le q-1$ for all $0<q<1$. Moreover, if we can show that $\tau(q)$ is strictly convex on $[0,1]$, this will imply that $\tau(q) < q-1$ for all $0<q<1$ and, as noted before (see (\ref{eq:dimdiff}) and the associated remark), this will yield the theorem.

A straightforward calculation shows that if we let
\[
\tau_k(q) = \frac{\log W_k(q)}{k\log(1/\lambda)},
\]
then
\begin{equation} \label{eq:derivtau}
\tau'_k(1) = \frac{\sum_{j=1}^{N_k} p_k(j) \log p_k(j)}{\log(\lambda^{-k})}.
\end{equation}
Denote the sum in the numerator above by $h_k$. In the course of the proof that $\dim_H(\nu_\lambda)<1$ (which essentially goes back to Garsia \cite{garsia}), it is shown that
\begin{equation} \label{eq:entropyineq}
\lim_{k\rightarrow\infty} \left(h_k - \log(N_k) \right)= -\infty.
\end{equation}
See \cite{prurb}, pp. 179-180 for a proof of this fact. We remark that the proof uses both Garsia's lemma and the singularity of $\nu_\lambda$, but is otherwise elementary.

Recalling that $N_k$ is bounded by a constant multiple of $\lambda^{-k}$, we deduce from (\ref{eq:derivtau}) and (\ref{eq:entropyineq}) that $\tau'_{k_0}(1) < 1$ for some sufficiently large $k_0$. Since, on the other hand, $\tau_{k_0}(1)=0$ and $\tau(q)\le \tau_{k_0}(q)$ for all $0\le q\le 1$, we conclude that $\tau'_{-}(1)<1$, where $\tau'_{-}(1)$ denotes the left derivative. This shows that $\tau(q)$ cannot agree with $1-q$ on $(0,1)$ and therefore must be strictly convex, completing the proof $\blacksquare$

We make some remarks about the above proof. First, the proof is about the Hausdorff dimension of $K_{\gamma,\lambda}$. Under strong separation, the Box dimension does not drop when $\lambda^{-1}$ is Pisot. When there are overlaps, it is no longer so clear what happens to the Box dimension, but in principle there is no reason to believe it will also drop.

Second, the function $\tau(q)$ that appeared in the course of the proof is actually (the negative of) the $L^q$-spectrum of $\nu_\lambda$. The fact that $\tau(q)$ is strictly convex corresponds, under the multifractal formalism, to the fact that $\nu_\lambda$ is a multifractal measure; i.e. it has a range of local dimensions. In general, the left and right derivatives of the $L^q$-spectrum at $1$ give substantial information about the measure; see for example \cite{diffdimensions}. Thus the proof is another indication of the delicate relationship between Bernoulli convolutions and the sets $K_{\gamma,\lambda}$.

\begin{corollary}
Fix $\lambda$ such that $\lambda^{-1}$ is a Pisot number.
\begin{enumerate}
\item
If $\lambda>1/2$ then
\[
\dim_H(K_{\gamma,\lambda}) < 1 + \frac{\log(2\lambda)}{\log(1/\gamma)}
\]
for all $0<\gamma<\min(\lambda,1/(2\lambda))$.
\item
If $\lambda>1/\sqrt{2}$ then there exists $\varepsilon$ such that
\[
\dim_H(K_{\gamma,\lambda}) < 2
\]
for all $\gamma \in ((2\lambda)^{-1},(2\lambda)^{-1}+\varepsilon)$.

\end{enumerate}
\end{corollary}
\textit{Proof}. The first part is clear from Theorem \ref{th:dimdrop}. For the second part, note that when $\gamma = (2\lambda)^{-1}$, we have that $\gamma<\lambda$ (since $\lambda>1/\sqrt{2}$) and the Falconer dimension is exactly $2$. Therefore the second part follows from the continuity of the drop in Theorem \ref{th:dimdrop}. $\blacksquare$

Interestingly, there are exactly two Pisot numbers whose reciprocals are greater than $1/\sqrt{2}$. The smallest Pisot number is the real root of $x^3-x-1$, which is about $1.324717\ldots$. The second smallest Pisot number is the positive root of the polynomial $x^4-x^3-1$; it is about $1.380280\ldots$. The next Pisot number is the positive root of $x^5-x^4-x^3+x^2-1$, which is already greater than $1.44327 > \sqrt{2}$. See \cite{pisotsalem}, Theorem 7.2.1. for a proof of these facts.

 \subsection{Zero Hausdorff measure}

In this subsection we assume $D=\{0,1,\ldots,m-1\}$. We showed that when the Falconer dimension is less than $2$, the Hausdorff dimension of $K_{\gamma,\lambda}$ is almost everywhere equal to the Falconer dimension. Hence it is natural to ask what is the Hausdorff measure in the critical dimension $s=1+\log(m\lambda)/\log(1/\gamma)$. Unlike the self-similar case, this is a non-trivial question even in the strong separation case, provided $\lambda>1/m$ (otherwise $K_{\gamma,\lambda}$ can be seen to be bi-Lipschitz equivalent to the Cantor set $K_\lambda$). It is in fact very easy to show that $\mh^s(K_{\gamma,\lambda})<\infty$ for \textit{all} $(\gamma,\lambda)$, by considering the natural cover. The following result was communicated to us by M. Rams, but seems to be folklore.

\begin{theorem} \label{th:hmnooverlaps}
Assume the strong separation is verified for $K_{\gamma,\lambda}$, where $\lambda>1/m,\gamma\lambda<1/m$, and let $s$ be the Falconer dimension of $K_{\gamma,\lambda}$. Then $\mh^s(K_{\gamma,\lambda})>0$ if and only if $\nu_\lambda$ is absolutely continuous with a bounded density.
\end{theorem}
For the proof of the theorem we need the following lemma which has independent interest.
\begin{lemma} \label{lemma:equivalence}
Assume that $\lambda>1/m$, $\gamma\lambda<1/m$ and $0< \mh^s(K_{\gamma,\lambda})<\infty$, where $s$ is the Falconer dimension of $K_{\gamma,\lambda}$.  Then $\wmh^s|_{K_{\gamma,\lambda}}$ assigns the same mass to all cylinders of level $k$, where $\wmh^s$ is the measure of Hausdorff type obtained by considering covers by open squares only. (Here we do not assume a separation condition).
\end{lemma}
\textit{Proof of lemma}. For notational convenience we will omit the subscripts $\gamma,\lambda$. Let $M=\wmh^s(K)$. We will argue by contradiction. Since
\[
\wmh^s(K) \le \sum_{u\in D^k} \wmh^s(K(u)),
\]
it follows that some $K(u)$ has measure \textit{larger} than $M m^{-k}$; fix such a $u$ and choose $\varepsilon>0$ such that
\begin{equation} \label{eq:lemmaequiv1}
\wmh^s(K(u))>(1+\varepsilon)M m^{-k}.
\end{equation}
By decomposing $u$ into sub-cylinders and using subadditivity again, it follows that for all sufficiently large $k$ there is $u\in D^k$ verifying (\ref{eq:lemmaequiv1}).

Now fix $\delta>0$ and choose a cover $\mathcal{C}=\{S_j\}$ of $K$ by squares such that
\[
\sum_j \diam(S_j)^s < (1+\delta) M,
\]
Also fix $k$ such that
\begin{equation} \label{eq:lemmaequiv2}
\gamma^{ks} = \left(\frac{\gamma}{m\lambda}\right)^k   < \varepsilon m^{-k}
\end{equation}
and (\ref{eq:lemmaequiv1}) holds for some $u\in D^k$. Consider a cover $\mathcal{C}'$ of $K(u)$ defined as follows: for each $j$, cover the rectangle $\phi_u(S_j)$ by $\lfloor \lambda^k/\gamma^k \rfloor+1$ squares of side $\gamma^k$, and take the union of those squares over all $j$. Therefore we have
\begin{equation} \label{eq:lemmaquiv3}
\wmh^s(K(u)) \le (\lfloor \lambda^k/\gamma^k \rfloor+1)\,\gamma^{ks} \sum_j \diam(S_j)^s < M (1+\delta)( m^{-k} + \gamma^{ks}),
\end{equation}
where we used that $\lambda\gamma^{1-s}=m^{-1}$. Letting $\delta\rightarrow 0$ and recalling (\ref{eq:lemmaequiv2}) we conclude that $\wmh^s(K(u))\le M(1+\varepsilon)m^{-k}$. This contradicts (\ref{eq:lemmaequiv1}), as desired. $\blacksquare$

\medskip

\textit{Proof of Theorem \ref{th:hmnooverlaps}}. We use the same notation as in the lemma. For $\omega\in\ D^\mathbb{N}$ let $S_k(\omega)$ be the open square centered at $\Pi(\omega)$ and half-side $\gamma^k/(1-\lambda)$. One consequence of strong separation, the self-affine relation and the fact that $\mu$ projects onto the Bernoulli convolution $\nu_\lambda$ is that
\begin{equation} \label{eq:hmnooverlaps1}
\mu(S_k(\Pi(\omega))) = 2^{-k} \nu_\lambda(B(\Pi_\lambda(\sigma^k\omega),(\gamma/\lambda)^k/(1-\lambda))).
\end{equation}
From this it immediately follows that if $\nu$ has a bounded density then
\[
\mu(S_k(\Pi(\omega))) < C (2\gamma/\lambda)^k = C \gamma^{ks},
\]
where $C$ is independent of $k$. Hence $\mu$ is an $s$-dimensional Frostman measure and it follows that $\mh^s(K)>0$.

Now suppose that $d\nu/dx$ is not a bounded function (or $\nu$ is not absolutely continuous at all), and fix $M>0$. Let
\[
E = \{ \omega\in D^\mathbb{N}: \nu_\lambda(B(\Pi_\lambda(\omega),r))>  M r \textrm{ for all } r<r_0  \}.
\]
If $r_0$ is small enough then $E$ has positive measure. But in this case, the ergodic theorem implies that for almost every $\omega$, $\sigma^k\omega$ visits $E$ with positive frequency, and therefore we deduce from (\ref{eq:hmnooverlaps1}) that
\[
\limsup_{k\rightarrow\infty} \frac{\mu(S_k(x))}{\gamma^{ks}} = \infty \quad \mu-a.e.\,x.
\]
Under strong separation, $\mu$ assigns the same mass to all cylinders of the same level; therefore $\mu$ is a constant multiple of $\wmh^s|_K$ by the lemma, and thus equivalent to $\mh^s|_K$. But by the density theorems (see \cite{mattila}, Theorem 6.2. (1)),
\[
\limsup_{k\rightarrow\infty} \frac{\mh^s(S_k(x)\cap K)}{\gamma^{ks}} \le C < \infty \quad \mh^s-a.e.\, x;
\]
therefore $\mh^s(K)=0$ as desired. (Note that we do need the lemma; otherwise $\mh^s$ might be concentrated in the exceptional set of $x$).  $\blacksquare$

\medskip

We remark that in the range $m^{-1}<\lambda<m^{-1/2}$ nothing is known about the boundedness of $d\nu_\lambda/dx$ generically, so we in fact do not know what the case is for $\mh^s(\mu_{\gamma,\lambda})$.

We now consider the overlapping case. In the self-similar setting, Solomyak \cite{measanddim} in a particular situation and later Peres, Simon and Solomyak \cite{zerohauspositpack} in greater generality showed that, assuming transversality, self-similar sets with overlap have typically $0$ measure. It turns out that the proof in \cite{zerohauspositpack} extends to our setting. Unfortunately, we are not able to check the needed concept of transversality (which is different from the one used in the first sections) in all of the relevant region $\Gamma(0,1/m)$, although it does hold for a large chunk of the overlapping region by results of Solomyak \cite{connectloci}.

The proof in \cite{zerohauspositpack} relies on the Bandt-Graf criterion \cite{bandtgraf}, so our first lemma extends one direction of this important criterion to our family of self-affine sets.

\begin{lemma} \label{lemma:bandtgraf}
Let us say that two cylinders $K(u),K(v)$ are $\varepsilon$-relatively close if
\begin{equation} \label{eq:epsilonclose}
\|(\phi_u)^{-1}\phi_v-I\|<\varepsilon,
\end{equation}
where $I$ is the identity map and $\|\cdot\|$ denotes Euclidean operator norm.

If for every $\varepsilon$ there are finite words $u,v$ such that $K(u)$ and $K(v)$ are $\varepsilon$-relatively close then $\mh^s(K)=0$, where $s$ is the Falconer dimension.
\end{lemma}
\textit{Proof}. Our proof follows the idea of Bandt and Graf; details are provided for completeness. Assume by way of contradiction that $\mh^s(K)>0$. We know from Lemma \ref{lemma:equivalence} that in this case all cylinders of the same level are disjoint in $\wmh^s$-measure.

Let $M=\wmh^s(K)$, and choose a cover $\{S_j\}$ of $K$ by open squares such that
\[
\sum_j \diam(S_j)^s < M +\delta,
\]
where $\delta$ is to be chosen later. Let $S$ be the union of the $S_j$, and pick $\varepsilon>0$ such that if $\|A-I\|<\varepsilon$ for some map $A$ then $AK\subset S$; this is possible since $S$ is open.

Next take $u,v$ such that (\ref{eq:epsilonclose}) holds with this $\varepsilon$. Note that if $\varepsilon$ is small enough then $u$ and $v$ must have the same length, say $k$; note also that $k$ can be made arbitrarily large by taking $\varepsilon$ small. Because of the way $\varepsilon$ was chosen, we have $K_u \subset \phi_v(S)$. Now adapt the covering $\{S_j\}$ to $K(v)$ as in Lemma \ref{lemma:equivalence} (mapping the $S_j$ by $\phi_v$ and then dividing the resulting rectangles into squares). The union of this covering clearly contains $\phi_v(S)$, whence it is also a covering of $K_u$; this is a contradiction since $K_u$ and $K_v$ are measure disjoint and this covering is almost optimal. More precisely, use (\ref{eq:lemmaquiv3}) to get
\[
2 M m^{-k} = \wmh^s(K_u) + \wmh^s(K_v) \le  M(1+\delta)(m^{-k}+\gamma^{sk}).
\]
Taking $\delta=1/2$ and $k$ large enough so that $\gamma^{ks} < (1/3) m^{-k}$ yields the desired contradiction. $\blacksquare$

In applying the previous lemma we will need to prove that certain sets have zero measure, and following \cite{zerohauspositpack} we will do so by showing that those sets have no Lebesgue density points. However, we will need a different notion of density points, defined using averages over rectangles of size $\gamma^k\times\lambda^k$ rather than balls. The following lemma shows that this makes no difference to us.

\begin{lemma} \label{lemma:densityrect}
Let $R_C((x,y),k)$ be the rectangle centered at $(x,y)\in (0,1)^2$ having dimensions $2C x^k\times 2C y^k$ ($C$ is allowed to depend on $(x,y)$ but not on $k$). Let $A$ be a measurable subset of $(0,1)^2$. If for all $u=(x,y)\in A$ we have that
\[
\limsup_{k\rightarrow\infty} \frac{\leb_2(R_C(u,k)\cap A)}{\leb_2(R_C(u,k))} < 1,
\]
then $\leb_2(A)=0$.
\end{lemma}
\textit{Proof}. This is a special case of the fact that the set of rectangles with sides parallel to the axes is a density basis of $\mathbb{R}^2$; see \cite{guzman}, Theorem 3.1. $\blacksquare$

Now we define the appropriate notion of transversality for this setting.

\begin{definition}
We say that $\mathcal{R}\subset\mathbb{R}^2$ is a region of $\star$-transversality if for all $(\gamma,\lambda)$ there is $f\in\mb$ such that $f(\gamma)=f(\lambda)=0$ but $f'(\gamma)\ne 0$, $f'(\lambda)\ne 0$.
\end{definition}
The requirement that $\gamma$ and $\lambda$ are zeros of $f$ is natural since we want $K_{\gamma,\lambda}$ to have overlaps (recall Lemma \ref{lemma:conneclocus}; also note that $K_{\gamma,\lambda}$ is not necessarily connected if $m>1$). Note that in the previous definition of transversality one of $\gamma$, $\lambda$ was allowed to be a double zero as long as the other parameter was not; this is not the case here. On the other hand, here we need only the \textit{existence} of such an $f$; it is therefore natural to conjecture that
\begin{equation}
\mathcal{R}=\Gamma(0,1/m)\cap \mathcal{IP}_m,
\end{equation}
where $\mathcal{IP}_m$ denote the set of parameters where there is an overlap ($\mathcal{IP}$ stands for ``Intersection Parameters''; this terminology was introduced in \cite{zerohauspositpack}). However, we were not able to take advantage of the fact that only existence of an $f$ is needed.

Solomyak \cite{connectloci} has obtained a large region of $\star$-transversality in the case $m=2$. Although it does not cover all of $\mathcal{IP}_2$, it does contain a big chunk of it. We refer to his paper for details.

\begin{theorem} \label{th:zeromeasure}
Let $\mathcal{R}\subset\mathbb{R}^2$ be a region of $\star$-transversality. Then for almost all $(\gamma,\lambda)\in \mathcal{R}$,
\[
\mh^{d(\gamma,\lambda)}(K_{\gamma,\lambda}) = 0,
\]
where $d(\gamma,\lambda)=1+\frac{\log (2\lambda)}{1/\log{\gamma}}$.
\end{theorem}
\textit{Proof}. The proof follows closely the proof of Theorem 2.1 in \cite{zerohauspositpack} (in fact, the easier homogeneous case). To begin, observe that for $\varepsilon>0$ small, two cylinders $K_{\gamma,\lambda}(u), K_{\gamma,\lambda}(v)$ are $\varepsilon$-relatively close if and only if $|u|=|v|=k$ for some large $k$ and
\begin{equation} \label{eq:zeromeas1}
|\Pi_{\alpha}(u)-\Pi_\alpha(v)| \le \varepsilon \alpha^k,\quad \alpha=\gamma,\lambda.
\end{equation}
Let $\Phi_\varepsilon$ be the set of all $(\gamma,\lambda)\in \mathcal{R}$ such that (\ref{eq:zeromeas1}) holds for some $u,v$ of the same length. Invoking Lemma \ref{lemma:bandtgraf}, it is enough to show that $\mathcal{R}\backslash \Phi_\varepsilon$ has zero Lebesgue measure for all $\varepsilon$ sufficiently small. To this end, we will show the following: for all $(\gamma_0,\lambda_0) \in \mathcal{R}$, some $C>0$, $0<\eta<1$, and all sufficiently large $k$ (all depending on $(\gamma_0,\lambda_0))$,
\begin{equation} \label{eq:conditionzeromeas}
\leb_2((\mathcal{R}\backslash \Phi_\varepsilon) \cap  R_C((\gamma_0,\lambda_0),k)) < (1-\eta) \leb_2 (R_C((\gamma_0,\lambda_0),k)),
\end{equation}
where $R_C((\gamma_0,\lambda_0),k))$ is as defined in Lemma \ref{lemma:densityrect}. Once we have shown this, the same lemma will give that $\mathcal{R}\backslash\Phi_\varepsilon$ has zero Lebesgue measure.

From now on fix $\varepsilon>0$ and $(\gamma_0,\lambda_0)\in \mathcal{R}$. Since $\mathcal{R}$ is a region of $\star$-transversality, there exists $\omega\in\mb$ such that $f(\gamma_0)=f(\lambda_0)=0, f'(\gamma_0)\ne 0, f'(\lambda_0)\ne 0$. Write $f=\omega_1 - \omega_2$, where $\omega_i \in D^\mathbb{N} (i=1,2)$. It follows that
\[
K(\omega_1|k) \cap K(\omega_2|k) \ne \varnothing.
\]
for all $k>0$ (where $\omega_i|k$ denotes the initial word of length $k$ of $\omega_i$). Let $U$ be a small open square centered at $(\gamma_0,\lambda_0)$ and compactly contained in $(0,1)^2$ such that
\begin{equation} \label{eq:c1}
c_1 := \frac{1}{2}\min\{ \min(|f'(\gamma)|,|f'(\lambda)|):(\gamma,\lambda) \in U \} > 0.
\end{equation}
Since the closure of $U$ is contained in $(0,1)^2$,
\begin{equation} \label{eq:c2}
c_2 := 2 \max\{ \max(|f'(\gamma)|,|f'(\lambda)|):(\gamma,\lambda) \in U \} < \infty.
\end{equation}
Also, since $f(\gamma_0)=f(\lambda_0)=0$,
\begin{equation} \label{eq:c3}
|\Pi_\alpha(\omega_1|k)-\Pi_\alpha(\omega_2|k)| \le c_3 \alpha^k\quad \alpha=\gamma_0,\lambda_0,
\end{equation}
for some finite $c_3$.

Let $g_k(x) = \Pi_x(\omega_1|k)-\Pi_x(\omega_2|k)$. Writing $g_k = f + (g_k-f)$ we see from (\ref{eq:c1}) and  (\ref{eq:c2}) that, if $k$ is large enough,
\begin{equation} \label{eq:derbounds}
(\gamma,\lambda)\in U\,\Longrightarrow\, c_1 < |g_k'(\gamma)|,|g_k'(\lambda)| < c_2,
\end{equation}

Since $U$ is open, (\ref{eq:c3}) and (\ref{eq:derbounds}) imply that, for large $k$, there exist $\gamma_1,\lambda_1$ such that $|\alpha_1-\alpha_0|< c_3 \alpha_0^k / c_1$ and $g_k(\alpha_1)=0$ for $\alpha=\gamma$ or $\lambda$.

For this choice of $\gamma_1,\lambda_1$, the upper bound in (\ref{eq:derbounds}) implies that
\[
0 < t < \varepsilon \alpha_1^k / c_2 \quad \Longrightarrow \quad |g_k(\alpha_1+t)|< \varepsilon \alpha_1^k < \varepsilon (\alpha_1+t)^k,\quad \alpha=\gamma,\lambda,
\]
whence
\begin{equation} \label{eq:T}
T := (\gamma_1,\gamma_1+\varepsilon\gamma_1^k / c_2) \times  (\lambda_1,\lambda_1+\varepsilon\lambda_1^k / c_2) \subset \Phi_\varepsilon.
\end{equation}
Note however that if $k$ is large enough
\[
\alpha_1^k > (\alpha_0^k - c_3 \alpha_0^k /c_1)^k = \alpha_0^k \left(1 - c_3 \alpha_0^{k-1}/c_1 \right)^j > (1-1/k)^k \alpha_0^k > \frac{\alpha_0^k}{2},
\]
where $\alpha=\gamma$ or $\lambda$. A similar argument shows that $\alpha_1^k < 2\alpha_0^k$ for $\alpha=\gamma$ or $\lambda$ and large enough $k$. Therefore
\[
\leb_2(T) =  c_2^{-2} \varepsilon^2 \gamma_1^k\lambda_1^k \ge 4^{-1}c_2^{-2} \varepsilon^2 \gamma_0^k\lambda_0^k = 16^{-1} c_2^{-2} C^{-2} \varepsilon^2 \leb_2(R_C((\gamma_0,\lambda_0),k)).
\]
On the other hand,
\[
|(\alpha_1 + t) - \alpha_0| \le  (c_3/c_1 + 2\varepsilon/c_2) \alpha_0^k \,\textrm{ for all } 0\le t\le \alpha_1^k/c_2 , \quad \alpha=\gamma,\lambda.
\]
Therefore if we take $C= c_3/c_1 + 2\varepsilon/c_2$ and $\eta=16^{-1} c_2^{-2} C^{-2} \varepsilon^2$, we obtain that
\[
T \subset R_C((\gamma_0,\lambda_0),k);\quad \leb_2(T) > \eta\,\leb_2(R_C((\gamma_0,\lambda_0),k)).
\]
Together with (\ref{eq:T}) this implies (\ref{eq:conditionzeromeas}), and the proof is complete. $\blacksquare$

\subsection{Laws of pure type} \label{sec:puretype}

It is known that self-similar measures are either singular or mutually absolutely continuous with respect to Lebesgue measure on the attractor. This was proved, in increasing levels of generality, in \cite{equivalencebernoulli}, \cite{sixtyyears} and \cite{abscontinuity}. All of those papers use the Lebesgue density theorem and assume that the maps are at least conformal; however, by using density bases more general than balls it is possible to adapt those proofs to our setting.

Recall the a density basis $\mathcal{V}$ for a Borel set $K$ is a family of open sets such that the following holds: for all $x\in K$ there are arbitrarily small sets $V\in\mathcal{V}$ containing $x$ and if $A\subset K$ is a Borel set, then
\[
\lim_{V\rightarrow x, x\in V\in\mathcal{V}} \frac{\leb_n(V\cap A)}{\leb_n(V)} = \mathbf{1}_A(x) \textrm{ for }  \leb_n-a.e. x\in \mathbb{R}^n.
\]
(Here $V\rightarrow x$ means that $\diam(V\cup\{x\})\rightarrow 0$).

Let $K$ be the attractor of an affine i.f.s. $\{\phi_1,\ldots,\phi_m\}$ on $\mathbb{R}^n$, and assume that $\leb_n(K)>0$. We will say that $K$ is \em differentiation-regular \em if there exists a density basis $\mathcal{V}$ for $K$ and a constant $\eta>0$ such that the following holds: for every $x\in K$ there is a sequence $\{V_j(x)\}$ in $\mathcal{V}$ with $x\in V_j(x), V_j(x)\rightarrow x$, such that if $V=V_j(x)$ for some $j$, there exists a finite word $u=u_j(x)$ verifying
\begin{equation} \label{eq:densitybasecondition}
\phi_{u}(K) \subset V;\quad \leb_n(\phi_{u}(K)) \ge \eta \leb_n(V).
\end{equation}

Very roughly speaking, a self-affine set is differentiation-regular if we can pick a differentiation basis for $K$ consisting of open sets which look like some cylinder in the construction of $K$. The next proposition shows that some important classes of self-affine sets, including those studied in this paper, are indeed differentiation-regular.

\begin{proposition} \label{prop:diffregular}
Let $\phi_1,\ldots,\phi_m$ be contracting affine maps on $\mathbb{R}^n$, and let $A_i$ denote the linear part of $\phi_i$. Let $K$ be the attractor of $\{\phi_1,\ldots,\phi_m\}$.
Assume that $\leb_n(K)>0$ and any one of the following conditions hold:
\begin{enumerate}
\item $n=2$ and all the maps $A_i$ are equal.
\item There exists a finite generating system $W$ of $\mathbb{R}^n$, such that $A_i W \subset W$ for all $1\le i\le m$.
\item All the maps $A_i$ are simultaneously diagonalizable.
\end{enumerate}
Then $K$ is differentiation-regular.
\end{proposition}
\textit{Proof}. Assume first that $n=2$ and all the linear parts are equal to $A$. Let $B$ be a ball centered at the origin and containing $K$, and let $\theta_k \in S^{1}$ be the direction of the major axis of the ellipse $A^k(B)$. Then we can pick a subsequence $\theta_{k_i}$ which is either constant or lacunary (a sequence $\{\theta_i\}\subset S^{1}$ is lacunary if it converges to some $l\in S^{n-1}$ and there is $C>1$ such that $C|\theta_{i+1}-l|<|\theta_i-l|$ for all $i$). In either case, it is well-known that the family
\[
\mathcal{V} = \{ A^{k_i}(B) + v : v\in\mathbb{R}^2, i\in\mathbb{N} \}
\]
is a density basis of $\mathbb{R}^2$ (This is originally due to R.Fromberg, a proof can be found in \cite{lacunary}). Now if $u$ is a word of length $k_i$ we have that
$\phi_u(K) \subset \phi_u(B) \in \mathcal{V}$ and
and
\begin{equation} \label{eq:diffregular}
\leb_2(\phi_u(K)) = \det(\phi_u)\leb_2(K) = \eta \det(\phi_u)\leb_2(B) = \eta \leb_2(\phi_u(B)),
\end{equation}
where $\eta=\leb_2(K)/\leb_2(B)$. This shows that $K$ is differentiation-regular.

Now we consider the second case. Fix a generating set $W$ as in the statement, and let $P$ be a convex polyhedra containing $K$ and whose sides are parallel to elements of $W$ (since $W$ contains a basis, we can take $P$ to be a suitable parallelepiped). The set of all convex polyhedra with sides parallel to some element of $W$ is known to be a density basis of $\mathbb{R}^n$ (\cite{guzman}, p.137); let us denote this basis by $\mathcal{V}$. By hypothesis, $\phi_u(K)\subset\phi_u(P)\in \mathcal{V}$ for all $u$, and we can conclude that $K$ is differentiation-regular as in (\ref{eq:diffregular}).

If all the $A_i$ are simultaneously diagonalizable then we can apply the previous case: just take $W$ to be the set of (simultaneous) eigenvectors. $\blacksquare$

Now we extend Proposition 3.1 in \cite{sixtyyears} to differentiation-regular self-affine sets.

\begin{proposition} \label{prop:puretype}
Let $K$ be the attractor of an affine i.f.s. $\{\phi_1,\ldots,\phi_m\}$ on $\mathbb{R}^n$, and assume that $K$ is differentiation-regular. Then for any self-affine measure $\mu$ supported on $K$, $\leb_n|K$ is either singular or mutually absolutely continuous with respect to $\mu$.
\end{proposition}
\textit{Proof}. It is well known that $\mu$ is either singular or absolutely continuous with respect to Lebesgue measure (this can be seen for instance by decomposing $\mu$ into absolutely continuous and singular parts, and showing that each of them also verifies the self-affine relation, hence one of them must be trivial). Therefore it is enough to show that if $\mu$ is the attractor of the weighted i.f.s. $\{ (\phi_i,p_i)\}_{i=1}^m$ and $\mu\ll\leb_n$, then $\leb_n|_K \ll \mu$.

Following \cite{sixtyyears}, Proposition 3.1, let
\[
\beta = \frac{1}{\leb_n(K)} \sup\{ \leb_n(A) : A \textrm{ Borel }, A\subset K, \mu(A)= 0\}.
\]
Since we are assuming that $\mu$ is not singular, $0\le \beta<1$; we want to show that in fact $\beta=0$. Fix a Borel subset $A_0$ of $K$ such that $\mu(A_0)=0$. Let $x\in K$ and $V=V_j(x)$ for some $x$. Pick a word $u$ such that (\ref{eq:densitybasecondition}) holds. The Borel set $\phi_u^{-1}(A_0\cap \phi_u(K))$ is contained in $K$; moreover it has zero $\mu$-measure since $\mu$ is self-affine (and thus $\mu\circ\phi_u^{-1}$ is dominated by a multiple of $\mu$).

By the definition of $\beta$ and the fact that $\phi_u$ is affine,
\begin{eqnarray}
\det(\phi_u)^{-1} \leb_n (A_0\cap \phi_u(K)) & =  & \leb_n(\phi_u^{-1}(A_0\cap \phi_u(K)) \nonumber\\
& \le &  \beta \leb_n(K) \nonumber\\
& = & \beta\det(\phi_u)^{-1} \leb_n(\phi_u(K))\nonumber.
\end{eqnarray}

Therefore, using (\ref{eq:densitybasecondition}),
\[
\leb_n(V\backslash A_0) \ge (1-\beta)\leb_n(\phi_u(K)) \ge (1-\beta)\eta \leb_n(V).
\]
The definition of density basis then implies that $\leb_n(A_0)=0$, as desired $\blacksquare$

We finish this section with some remarks. First, we did not really use that $\mu$ is self-affine; just that
\begin{equation} \label{eq:weaksa}
\mu(A) = 0 \Rightarrow \mu(\phi_i^{-1}(A)) = 0,
\end{equation}
for all Borel sets $A$. We may call measures verifying (\ref{eq:weaksa}) \em weakly self-affine \em.

Second, although one can check the differentiation-regular condition in many interesting cases, there are many natural instances where this property appears to fail (although proving it rigorously looks difficult). For instance, let $A_0,A_1$ be two closed angular sectors in the open first quadrant, disjoint except at the origin (for simplicity we consider only the case $n=2$). If the linear part of $\phi_i$ maps $A_i$ into $A_i$, then the i.f.s. $\{\phi_1,\phi_2\}$ induces another i.f.s. in the circle $S^1$, whose attractor is a Cantor set of directions. However, it is known that for many Cantor sets $C$, the set of rectangles with sides parallel to $C$ is \textit{not} a density basis \cite{cantordirections}, while no Cantor set is known for which the opposite is true. Thus checking out differentiation-regularity appears unlikely.

Of course, it could be that Proposition \ref{prop:puretype} is true for arbitrary self-affine sets, using a different method of proof. However we were not able to verify this and we believe that there may be a counterexample.

 \subsection{Non-collinear digits}

So far we have only considered self-affine sets where the digits lie all on the line $\{x=y\}$. When $D=\{0,1\}$ this is of course not a restriction, but when $m>1$ it is certainly a strong assumption. In fact, almost all the results in this paper can be generalized, in theory, to a more general setting described below; the problem is that the needed notion of transversality can be very hard to check, or simply false. On the plus side, most of the results can be shown to hold under a perturbation of the digits. In particular, transversality (in all of its various forms) is an open condition. If anything, this shows that collinearity of the digits is not a necessary condition for the type of results obtained in this paper.

From now and until the end of this subsection let $D_i=\{d_0^i,d_1^i,\ldots,d_{m-1}^i\}$ for $i=1,2$ two sets of digits. Let $K_{\gamma,\lambda}$ be the attractor of the i.f.s.
\[
\{\, \phi_i(x,y) =  (\gamma x+d_i^1,\lambda y+ d_i^2)\, \}.
\]
Let us write $\Omega_m=\{0,\ldots,m-1\}^\mathbb{N}$.

\begin{definition}
For $\omega\in \Omega_m$ let
\[
f_\omega^i(x) = \sum_{j=1}^\infty  d_{\omega(j)}^i x^j \quad i=1,2.
\]
We say that $R\subset\mathbb{R}^2$ is a \em region of transversality \em if whenever $(\gamma,\lambda)\in R$ and $\omega,\omega'\in \Omega_m$, either $\gamma$ is not a double root of $f_\omega^1$ or $\lambda$ is not a double root of $f_\omega^2$.
\end{definition}

In this framework, suitable versions of Theorems \ref{th:generalth} and \ref{th:generalthbis} apply. In general, obtaining regions of transversality looks very difficult unless the digits are almost collinear or have some very special form. However, the regions $\ms_2$ in the aforementioned theorems can still be efficiently estimated, since $I$ and $J$ are in this case independent from one another. As an example, we have the following result, whose proof works exactly as in Theorem \ref{th:generalth}.

\begin{theorem}
Let $J_i$ be open intervals such that $f^i_\omega$ has no double roots in $I_i$ for all $\omega\in \Omega_m; i=1,2$. Let $I_i$ be open intervals such that $\nu_\alpha^i$ has correlation dimension $1$ for almost every $\alpha\in I_i$, $i=1,2$ (here $\nu_\lambda^i$ is the B.C. associated to $D_i$). Finally, let $R$ be a region of transversality. Let
\begin{eqnarray}
\ms_1 & = & (I_1\times I_2) \cap R \cap \Gamma(0,1/m);\nonumber\\
\ms_2 & = & ((I_1\times J_2) \cup (I_2\times J_1)) \cap \Gamma(0,1/m);\nonumber\\
\ms & = & \ms_1\cup\ms_2\nonumber.
\end{eqnarray}
Then $\dim_H(K_{\gamma,\lambda}) = 1+\log(m\lambda)/\log(1/\gamma)$ for almost all $(\gamma,\lambda)\in\ms$. $\blacksquare$
\end{theorem}
\begin{corollary}
If $D_1,D_2$ are sufficiently small perturbations of $\{0,\ldots,m-1\}$ then $\dim_H(K_{\gamma,\lambda})$ has the expected value for almost every $(\gamma,\lambda)\in \Gamma(0,1/m)$. $\blacksquare$
\end{corollary}

Note that we are somehow in the reverse situation with respect to Falconer's theorem, which holds for \textit{ all } families of linear maps satisfying certain conditions and for \textit{ almost every } digit. Here, the result holds for \textit{ all } digits in some open set and \textit{ almost every } parameter.

Another possible variant is to allow the digits to depend on the parameters. This is standard in the self-similar setting and present no additional complications here; most of the results, when stated appropriately, hold also in this case.

\textbf{Acknowledgements}. The author wishes to thank Boris Solomyak for his guidance and help in the completion of this paper.

\newcommand{\etalchar}[1]{$^{#1}$}

\end{document}